\documentclass[oneside,11pt]{article}

\usepackage{amsmath, amsfonts, amssymb, amsthm}
\usepackage{bbm}                  

\usepackage{graphicx}             
\usepackage{float}                
\usepackage{xcolor}               

\usepackage[labelformat=simple]{subcaption} 
\DeclareCaptionLabelFormat{subcaptionlabel}{\normalfont(\textbf{#2})}
\usepackage{ifthen}
\newboolean{figcolor}             

\usepackage{algorithm}
\usepackage{algpseudocode}

\usepackage{enumitem}             

\usepackage[top=1.5cm,bottom=2cm,right=2cm,left=2cm]{geometry}
\usepackage{titling}              

\usepackage{natbib}               

\usepackage{xparse}               

\allowdisplaybreaks              

\newtheorem{proposition}{Proposition}
\newtheorem{lemma}{Lemma}
\newtheorem{corollary}{Corollary}
\newtheorem{remark}{Remark}

\newtheorem{assumption}{Assumption}
\renewenvironment{proof}{\textit{Proof.}}{\qed\\}

\newcommand{\R}{\mathbb{R}}       
\newcommand{\N}{\mathbb{N}}       
\newcommand{\cR}{\mathcal{R}}     
\newcommand{\cF}{\mathcal{F}}     
\newcommand{\cD}{\mathcal{D}}     
\newcommand{\cC}{\mathcal{C}}     
\newcommand{\cB}{\mathcal{B}}     
\newcommand{\cA}{\mathcal{A}}     
\newcommand{\F}{\mathbb{F}}       
\newcommand{\sT}{\mathsf{T}}      

\newcommand{\lp}{\left(}   \newcommand{\rp}{\right)}
\newcommand{\lb}{\left[}   \newcommand{\rb}{\right]}
\newcommand{\lc}{\left\{}  \newcommand{\rc}{\right\}}
  
\newcommand{\lrav}[1]{\left\lvert#1\right\rvert}
\newcommand{\lrp}[1]{\lp#1\rp}

\newcommand{\lrc}[1]{\lc#1\rc}

\newcommand*{\defeq}{\mathrel{\vcenter{
                     \baselineskip0.5ex\lineskiplimit0pt
                     \hbox{\scriptsize.}\hbox{\scriptsize.}
                     }}%
                     =}
\RenewDocumentCommand{\Pr}{oe{_^}}{
  \operatorname{\mathbb{P}}%
  \IfValueT{#2}{_{#2}}%
  ^{\IfValueTF{#3}{#3}{}}%
  \IfValueT{#1}{\lp #1 \rp}%
}
\NewDocumentCommand{\Esp}{oe{_^}}{%
  \operatorname{\mathbb{E}}%
  \IfValueT{#2}{_{#2}}%
  ^{\IfValueTF{#3}{#3}{}}%
  \IfValueT{#1}{\lb #1 \rb}%
}
\NewDocumentCommand{\Var}{oe{_^}}{%
  \operatorname{\mathbb{V}\mathrm{ar}}%
  \IfValueT{#2}{_{#2}}%
  ^{\IfValueTF{#3}{#3}{}}%
  \IfValueT{#1}{\lb #1 \rb}%
}
\NewDocumentCommand{\Cov}{oe{_^}}{%
  \operatorname{\mathbb{C}\mathrm{ov}}%
  \IfValueT{#2}{_{#2}}%
  ^{\IfValueTF{#3}{#3}{}}%
  \IfValueT{#1}{\lb #1 \rb}%
}
\newcommand{\InfGen}{\mathbb{L}}  
\newcommand{\Ind}{\mathbbm{1}}    
\newcommand{\supp}{\mathrm{supp}}  
\newcommand{\leqlr}{\leq_{\mathrm{lr}}} 

\newcommand{\rmd}{\mathrm{d}}      

\newcommand{\wtilde}[1]{\widetilde{#1}}
\newcommand{\what}[1]{\widehat{#1}}
\newcommand{\uline}[1]{\underline{#1}}
\newcommand{\oline}[1]{\overline{#1}}

\usepackage[pdftex,final,bookmarksnumbered,breaklinks,colorlinks]{hyperref}
\hypersetup{
	final,
	bookmarksnumbered=true,
	bookmarksopen=true,
	bookmarksopenlevel=0,
	unicode=false,
	pdftoolbar=true,
	pdfmenubar=true,
	pdffitwindow=false,
	pdftitle={On the optimal stopping of randomized Gauss--Markov bridges},
	pdfdisplaydoctitle=true,
	pdftoolbar=true,
	pdfmenubar=true,
	pdflang={English},
	pdfauthor={Abel Azze, Bernardo D'Auria},
	pdfsubject={arXiv paper},
	pdfcreator={Abel Azze},
	pdfproducer={Abel Azze},
	pdfkeywords={Gauss--Markov}{Free-boundary problem}{Optimal stopping},
	pdfnewwindow=true,
	breaklinks=true,
	hidelinks,
	linkcolor=black,
	citecolor=black,
	filecolor=magenta,
	urlcolor=cyan
}

\setboolean{figcolor}{true}

\begin{document}

\title{On the optimal stopping of randomized Gauss--Markov bridges}
\setlength{\droptitle}{-1cm}
\predate{}%
\postdate{}%
\date{}

\author{Abel Azze$^{1, 3}$ and Bernardo D'Auria$^{2}$}
\footnotetext[1]{Department of Quantitative Methods, CUNEF Universidad (Spain).}
\footnotetext[2]{Department of Mathematics ``Tullio Levi Civita'', University of Padova (Italy).}
\footnotetext[3]{Corresponding author. e-mail: \href{mailto:abel.guada@cunef.edu}{abel.guada@cunef.edu}.}
\maketitle

\begin{abstract}
    We consider the optimal stopping problem for a Gauss–Markov process conditioned to adopt a prescribed terminal distribution. By applying a time–space transformation, we show it is equivalent to stopping a Brownian bridge pinned at a random endpoint with a time‐dependent payoff. We prove that the optimal rule is the first entry into the stopping region, and establish that the value function is Lipschitz continuous on compacts via a coupling of terminal pinning points across different initial conditions. A comparison theorems then order value functions according to likelihood–ratio ordering of terminal densities, and when these densities have bounded support, we bound the optimal boundary by that of a Gauss–Markov bridge. Although the stopping boundary need not be the graph of a function in general, we provide sufficient conditions under which this property holds, and identify strongly log‑concave terminal densities that guarantee this structure. Numerical experiments illustrate representative boundary shapes.
\end{abstract}
\begin{flushleft}
	\small\textbf{Keywords:} Optimal stopping; diffusion process; Time-inhomogeneity, Brownian bridge, Gauss--Markov processes.
\end{flushleft}

\section{Introduction}

Optimal Stopping Problems (OSPs) have long played a central role in probability theory and its applications, ranging from financial option pricing to sequential analysis and decision-making under uncertainty. 
One key component of these types of problems is the underlying process that one seeks to optimally stop to maximize a given expected payoff.
Over the years, researchers have studied a variety of models, from the simple Brownian Motion (BM) and the traditional Black–Scholes framework to more elaborated dynamics such as the Heston stochastic‑volatility model. 
However, little has been said to capture the fundamental idea that, no matter how sophisticated a model may be, it is still a belief, and thus, it should remain flexible to adapt when new information becomes available. 

One natural form of additional information is the value of the process at a future date. For example, an investor might believe that an asset price $X_T$ at time $T$ follows a specified distribution and wish to incorporate this belief when choosing the optimal selling time. More precisely, let $X=(X_t)_{t\in[0,T]}$ be a stochastic process and let \(Z\) be a random variable with a given law. What are the dynamics of $X$ given $X_T=Z$, and how does this conditioning alter the associated optimal stopping problem?

The literature on these types of bridge‑conditioned OSPs is scarce and lacks generality. \cite{Shepp-1969-explicit} provided the first explicit solution by showing that the OSP for a Brownian Bridge (BB) can be transformed into an infinite‑horizon OSP of a BM via a time–space change of variables. Subsequent work has extended his result in several directions, such as proposing different solution methods (\cite{Ekstrom-2009-optimal}, \cite{Ernst-2015-revisiting}), using a non-identity gain function (\cite{Ekstrom-2009-optimal}, \cite{DeAngelis-2020-optimal}), and including discounts (\cite{DAuria-2020-discounted}). More recently, the work of \cite{Azze-2024-optimal-TDOU} and \cite{Azze-2024-optimal-GMB} generalized Shepp's transformation method to the class of Gauss–Markov Bridges (GMB). Other authors included randomization in the terminal time (\cite{Follmer_1972_optimal}, \cite{Glover-2022-optimally}), as well as in the process's value at the horizon (\cite{Leung_2018_optimal}, \cite{Ekstrom-2020-optimal}).

In this more general setting, where the terminal value of the Brownian bridge is prescribed by an arbitrary distribution, the resulting (information) drift becomes highly nonlinear in the spatial variable, time‐inhomogeneous, and may violate the usual global Lipschitz condition. This extra layer of complexity has hampered the derivation of general results. 
The work of \cite{Leung_2018_optimal} tackled the problem mainly numerically, showing that the optimal stopping boundary need not be the graph of a single function when the terminal law is two‐point discrete or double‐exponential continuous. \cite{Ekstrom-2020-optimal} later proved this disconnectedness in the two‐point case and obtained several general properties of the value function under arbitrary terminal distributions.  

In this paper, we study the OSP of a Gauss–Markov Process (GMP) conditioned to adopt a prescribed terminal distribution. 
Through a suitable time–space transformation, we demonstrate that the problem is equivalent to an OSP of a BB with a random pinning point, thereby extending the work of \cite{Leung_2018_optimal} and \cite{Ekstrom-2020-optimal}.
We show that the optimal stopping time is given by the first entry into the stopping region, and prove that the associated value function is Lipschitz continuous on compact sets. To establish this regularity, we introduce a novel construction that embeds processes, starting from different initial conditions, into a common probability space equipped with a convenient copula distribution over their terminal values.

We derive comparison results showing that the value functions are ordered when the corresponding terminal densities hold the same ordering in the likelihood–ratio sense. In particular, where the terminal density has bounded support, the Optimal Stopping Boundary (OSB) can be bounded by that of a GMB, comprehensively studied by \cite{Azze-2024-optimal-GMB}. 

While the OSB in this general setting may fail to be the graph of a single function, we provide sufficient conditions under which this property holds. In particular, we identify a class of strongly log-concave terminal densities for which the boundary exhibits this structure. Numerical experiments are presented to illustrate the shape of the OSB under representative parameter configurations, and all code is made available for reproducibility.

The remainder of the paper is organized as follows. In Section~\ref{sec:rGMB}, we introduce the notion of randomized GMBs and review key properties of GMPs. Section~\ref{sec:OSP} establishes the equivalence between the original OSP and a transformed problem involving a BB with a random pinning point. Section~\ref{sec:OSP-analysis} focuses on the analysis of the transformed problem. In Section~\ref{sec:single_boundary}, we provide sufficient conditions for the OSB to be the graph of a single function and solve the problem explicitly for Gaussian and degenerate terminal laws. Finally, Section~\ref{sec:numerics} presents numerical simulations that illustrate our theoretical results and offer insights into the geometry of the stopping and continuation regions.

\section{Randomized Gauss--Markov bridges}\label{sec:rGMB}

\subsection{Gauss--Markov processes}

Consider a standard BM $B = (B_t)_{t\in\R_+}$ defined in the filtered space $(\what{\Omega}, \what{\cF}, (\what{\cF}_t)_{t\geq0}, \Pr)$, where $\Pr$ is the Brownian measure and $(\what{\cF}_t)_{t\geq 0}$ is the natural filtration of $B$. For a terminal time $T > 0$ and an initial condition $x_0\in\R$, let $X = (X_t)_{t\in[0, T]}$ be the unique strong solution of the Stochastic Differential Equation (SDE)
\begin{align}\label{eq:SDE_GM}
    \left\{
    \begin{aligned}
        \rmd X_t &= \lrp{\alpha(t) + \beta(t)X_t}\rmd t + \zeta(t)\rmd B_t, \\
        X_0 &= x_0,
    \end{aligned}
    \right.
\end{align}
where $\alpha:[0, T]\rightarrow \R$, $\beta:[0, T]\rightarrow \R$, and $\zeta:[0, T]\rightarrow (0, \infty)$ are continuous functions. Such a process $X$, besides being Markovian, is also Gaussian, with a normal transition probability density (see, e.g., \cite{Buonocore-2013-some}) defined by the mean and variance functions
\begin{align}
    m_{t, x}(t') &\defeq \Esp[X_{t'} | X_t = x] =  \varphi(t, t')\lrp{x + \int_t^{t'} \frac{\alpha(u)}{\varphi(t, u)}\, \rmd u}, \label{eq:mean_GMP} \\
    v_t^2(t') &\defeq \Var[X_{t'} | X_t = x] =  \varphi^2(t, t') \int_t^{t'} \frac{\zeta^2(u)}{\varphi^2(t, u)} \, \rmd u \label{eq:var_GMP},
\end{align}
with $0\leq t < t' \leq T$, and for
\begin{align}
    \varphi(t, t') \defeq \exp\lrc{\int_t^{t'} \beta(u) \, \rmd u}. \label{eq:varphi_GMP}
\end{align}

It is also known (see, e.g., \cite{Borisov_1983_criterion}) that the (auto)covariance of $X$ admits the ``factorization''
\begin{align*}
    \Cov[X_t,X_t'] = r_1(\min(t, t'))r_2(\max(t, t')),
\end{align*} 
for (see, e.g., \cite{Buonocore-2013-some}) the unique-up-to-multiplicative-constants functions
\begin{align*}
    r_1(t) = \varphi(0,t)\int_0^t\frac{\zeta^2(u)}{\varphi^2(0,u)}\,\rmd u,\quad r_2(t) = \varphi(0,t).
\end{align*}
Moreover (see, e.g., \cite{Mehr-1965-certain}), $X$ admits the BM representation
\begin{align}\label{eq:BM_representation}
    X_t &= m_{0,x_0}(t) + r_2(t)B_{h(t)} 
    = m(t) + \varphi(t)\lrp{B_{h(t)} + x_0},
\end{align}
for $h(t) \defeq r_1(t)/r_2(t)$, $m(t) \defeq m_{0, 0}(t)$, and $\varphi(t)\defeq \varphi(0, t)$.

\subsection{The standard BM representation}

Note that $h(t) = r_1(t)/r_2(t)$ is a strictly increasing function ranging from $0 = h(0)$ to $\mathsf{T} \defeq h(T)$. 
Hence, by defining the change-of-time and auxiliary functions
\begin{align}\label{eq:change-time}
    s = s(t) = h(t)/\sT, \quad 
    a_0(s) \defeq m(h^{-1}(s)), \quad 
    a_1(s) \defeq \varphi( h^{-1}(s))\sqrt{\sT},
\end{align}
we obtain from \eqref{eq:BM_representation} and the scale property of a BM, that
\begin{align}\label{eq:BM_representation_process}
    X_t = G(s, Y_s), \quad s\in[0, 1],
\end{align}
where $(Y_s)_{s\in[0, 1]}$ is a BM starting at $Y_0 = x_0$, and $G$ takes the form
\begin{align}\label{eq:gain_function}
    G(s, y) \defeq a_0(s) + a_1(s)y.
\end{align}
From the definition of $a_0$ and $a_1$ in \eqref{eq:change-time} and \eqref{eq:mean_GMP}-\eqref{eq:varphi_GMP}, and the continuity of $\alpha$, $\beta$ and $\zeta$ in $[0, T]$, it follows
that
\begin{align*}
    a_0 \text{ and } a_1 \text{ are continuously differentiable in } [0, 1], \\
    a_1(s) > 0 \text{ for all } s \in [0, 1],
\end{align*}
and the following positive constants are well defined:
\begin{align}\label{eq:bounds}
    A_0 &\defeq \sup_{s\in\R_+} \lrav{a_0(s)},\quad A_0' \defeq \sup_{s\in\R_+} \lrav{a_0'(s)},\quad 
    A_1 \defeq \sup_{s\in\R_+} a_1(s),\quad A_1' \defeq \sup_{s\in\R_+} \lrav{a_1'(s)}.
\end{align}

\subsection{Conditioning on a terminal density}\label{sec:terminal_density}

For a prescribed density function $\nu$ in $\R$, consider a random variable $Z$ living in the probability space $(\wtilde{\Omega}, \wtilde{\cF}, \wtilde{\Pr}^\nu)$, where $\wtilde{\Pr}^\nu$ is the probability measure associated to $\nu$.
Construct the filtered space $(\Omega, \cF, \cF_t, \Pr^\nu)$ such that $\Omega = \what{\Omega}\times\wtilde{\Omega}$, $\cF = \what{\cF}\otimes\wtilde{\cF}$, $\cF_t = \what{\cF}_t\vee\sigma(Z)$, and $\Pr^\nu(\cdot) = \Pr\otimes\wtilde{\Pr}^\nu(\cdot \mid X_T = G(1, Z)) = \what{\Pr}\otimes\wtilde{\Pr}^\nu(\cdot \mid Y_1 = Z)$.
Hence, the process $(Y_{s})_{s\in[0, 1]}$, living in the filtered space $(\Omega, \cF,\F,\Pr^\nu)$, is such that $Y_1$ is distributed according to $\nu$, while $X_T$ is distributed like $\wtilde{\nu}(z) \defeq \nu((z - a_0(1))/a_1(1))/a_1(1)$. 
By means of the Girsanov theorem (see Proposition \ref{pr:bridge_SDE}), the process $(Y_{s})_{s\in[0, 1]}$ solves the SDE
\begin{align}\label{eq:SDE_rBB}
    \left\{
    \begin{aligned}
        \rmd Y_s &= \mu^\nu(s, Y_s) \rmd s + \rmd B_s^\nu, & \\
        Y_0 &= x_0,
    \end{aligned}
    \right.
\end{align}
for the $\Pr^\nu-$standard BM $B^\nu = (B_s^\nu)_{s\in[0,1]}$, and for the drift function 
\begin{align}\label{eq:drift_rBB}
    \mu^\nu(s,y) = \partial_y\ln(\psi^\nu(s, y)),
\end{align}
where $\psi^{\nu}$ is the Radon–Nikodym density $\Pr^\nu/\Pr$, defined as
\begin{align}\label{eq:RN_derivative}
    \psi^{\nu}(s, y) = \psi^{\nu}(s, y \mid 0, x_0) \defeq  
    \int_\R \frac{\phi(z;y,1-s)}{\phi(z;x_0, 1)}\nu(z)\,\rmd z,
\end{align}
where $\phi(z;\theta,\gamma^2)$ represents the normal density with mean $\theta$ and variance $\gamma^2$. 
Here and thereafter, we use $\partial_s$ and $\partial_y$ to refer to the partial derivatives to respect to the time and the space components, and $\partial_{yy}$ to indicate the second partial spacial derivative. 

For a comprehensive treatment of SDEs obtained by conditioning processes to adopt a prescribed terminal distribution, see \cite{Baudoin-2002-conditioned} and \cite{Macrina-2021-stochastic}.

\subsection{Randomized Brownian bridge}

From the random nature of the pinning point, we call the process defined at \eqref{eq:SDE_rBB} a randomized Brownian Bridge (rBB). Analogously, the more general class of processes $X$ defined in \eqref{eq:SDE_GM} and considered under the probability measure $\Pr^\nu$ are called randomized Gauss--Markov Bridges (rGMBs). 
A straightforward differentiation yields
\begin{align}\label{eq:SDE_rBB_non-Markov}
    \mu^\nu(s, y) = \partial_y \ln\lrp{\psi^{\nu}(s, y)}
    &= \frac{
    \int_\R \frac{z-y}{1-s} \frac{\phi(z;y,1-s)}{\phi\lrp{z;x_0,1}}\nu(z)\,\rmd z
    }{
    \int_\R \frac{\phi(z;y,1-s)}{\phi(z;x_0,1)}\nu(z)\,\rmd z
    } 
    = \frac{\Esp[Z_{s,y}]^\nu - y}{1-s},
\end{align}
where $\Esp^\nu$ is the expectation with respect to $\Pr^\nu$, and $Z_{s,y}$ is the random pinning point provided $Y_s = y$, that is, $\Pr[Z_{s,y} \in \cdot]^\nu = \Pr[Z \in \cdot \mid Y_s = y]^\nu$. 
Note that the density of $Z_{s,y}$ under $\Pr^\nu$ is given by 
\begin{align}\label{eq:updated_terminal_dens}
    \nu_{s,y}(z) = \partial_z \Pr[Z \leq z | Y_s = y]^\nu = \frac{\phi(z;y,1-s)}{\phi(z;x_0,1)}\nu(z) \Big/ \psi^{\nu}(s,y).
\end{align}
    
In the sequel, we rely on the following assumption.
\begin{assumption}[Regularity of the terminal density]\label{asm:terminal_density} \ \\
    For all $\varepsilon \in (0, 1)$, there exists a constant $L_\varepsilon > 0$ such that $\Var[Z_{s,y}] \leq L_\varepsilon$ for all $(s,y)\in(0,1-\varepsilon]\times \R$.
\end{assumption}

The next proposition obtains a Lipschitz-type bound for the Wasserstein-$1$ distance of the pinning points for different initial conditions. This property will prove itself useful to obtain bound of the partial derivatives of the value function  associated to a rBB.

\begin{proposition}[Wasserstein-$1$-Lipschitz continuity of the pinning point] \label{pr:wasser-distance}\ \\
    For $s_1$, $s_2\in[0, 1)$ and $y_1$, $y_2 \in\R$, let $\pi^* = \pi_{(s_1,y_1),(s_2,y_2)}^*$ be Wasserstein-$1$ copula of $\nu_{s_1, y_1}$ and $\nu_{s_2, y_2}$. That is,
    \begin{align}\label{eq:wasser-1-distance}
        \Esp[|Z_{s_1,y_1} - Z_{s_2,y_2}|]^{\pi^*}
        = \mathcal{W}(\nu_{s_1,y_1}, \nu_{s_2,y_2}) 
        \defeq \inf_{\pi\in\Pi(\nu_{s,y_1}, \nu_{s,y_2})}\Esp[|Z_{s,y_1} - Z_{s,y_2}|]^{\pi},
    \end{align}
    where $\Pi(f, g)$ is the set of all couplings whose marginal densities are $f$ and $g$. Then, under Assumption \ref{asm:terminal_density}, 
    \begin{enumerate}[label=(\textit{\roman{*}}), ref=(\textit{\roman{*}})]
        \item for all $\varepsilon\in(0, 1)$, there exists a constant $L_\varepsilon > 0$ such that
        \label{pr:space-Lipschitz_pinning_point}
        \begin{align}\label{eq:space-Lipschitz_pinning_point}
            \Esp[|Z_{s, y_1} - Z_{s, y_2}|]^{\pi^*} \leq L_\varepsilon|y_1 - y_2|,
        \end{align} 
        for all $s\in[0, 1-\varepsilon]$ and all $y_1,y_2\in\R$, with $\pi^* = \pi_{(s,y_1),(s,y_2)}^*$.

        \item for all compact sets $\cR \subset [0, 1)\times \R$, there exists a constant $L_\cR > 0$ such that
        \label{pr:time-Lipschitz_pinning_point}
        \begin{align}\label{eq:time-Lipschitz_pinning_point}
            \Esp[|Z_{s_1, y} - Z_{s_2, y}|]^{\pi^*} \leq L_\cR|s_1 - s_2|,
        \end{align} 
        whenever $(s_1, y)$, $(s_2, y)\in\cR$, with $\pi^* = \pi_{(s_1,y),(s_2,y)}^*$.
    \end{enumerate}
\end{proposition}

\begin{proof}
    For both parts of the proposition, we rely on the following standard bound (see Proposition 7.10 from \cite{Villani-2003-topics}) of the Wasserstein-$1$ distance: 
    \begin{align}\label{eq:wasser-1_bound}
        \mathcal{W}(f, g) \leq \int_{\R} |z||f(z) - g(z)|\,\rmd z.
    \end{align}
    
    \ref{pr:space-Lipschitz_pinning_point}
    Straightforward differentiation of \eqref{eq:updated_terminal_dens} yields that 
    $\partial_y\nu_{s,y}(z) = \nu_{s,y}(z)\frac{z - \Esp[Z_{s, y}]^\nu}{1-s}$, 
    from where it follows that
    \begin{align*}
        \int_{\R}|z||\partial_y\nu_{s,y}(z)|\,\rmd z 
        &\leq \frac{1}{1-s}\int_{\R}\nu_{s,y}(z)|z||z - \Esp[Z_{s, y}]^\nu|\,\rmd z \\
        &\leq \frac{1}{1-s}\lrp{\int_{\R}\nu_{s,y}(z)(z - \Esp[Z_{s, y}]^\nu)^2\,\rmd z 
          + |\Esp[Z_{s, y}]^\nu|\int_{\R}\nu_{s,y}(z)|z - \Esp[Z_{s, y}]^\nu|\,\rmd z} \\
        &\leq \frac{1}{1-s}\lrp{\Var[Z_{s,y}]^\nu + |\Esp[Z_{s, y}]^\nu|\sqrt{\Var[Z_{s, y}]^\nu}},
    \end{align*}
    where, in the second inequality, we used the triangular inequality $|z| \leq |z - \Esp[Z_{s, y}]| + |\Esp[Z_{s, y}]|$.  Hence, the proof is an immediate consequence of Assumption \ref{asm:terminal_density}, inequality \eqref{eq:wasser-1_bound}, and the fact that
    \begin{align*}
        \int_{\R} |z||\nu_{s,y_1}(z) - \nu_{s,y_2}(z)|\,\rmd z 
        &\leq \int_{\R} |z|\int_{y_1}^{y_2}|\partial_y\nu_{s,y}(z)|\,\rmd y\,\rmd z 
        \leq \int_{y_1}^{y_2} \int_{\R}|z||\partial_y\nu_{s,y}(z)|\,\rmd z\,\rmd y.
    \end{align*}

    \ref{pr:time-Lipschitz_pinning_point} After differentiating \eqref{eq:updated_terminal_dens} with respect to time and performing some algebraic manipulations we obtain that
    \begin{align*}
        \partial_s\nu_{s,y}(z) = \nu_{s,y}(z)\frac{\Esp[(Z_{s,y} - y)^2]^\nu - (z-y)^2}{2(1-s)^2},
    \end{align*}
    from where it follows that there exists a bound $L_\cR > 0$, uniform for all $(s, y)\in \cR$, and controlled by the moments of $|Z_{s,y}|$ up to the third order, such that
    \begin{align*}
        \int_{\R}|z||\partial_y\nu_{s,y}(z)|\,\rmd z \leq L_\cR.
    \end{align*}
    Notice that here there is no need of using Assumption \ref{asm:terminal_density} since the existence and continuity in $\cR$ of the functions $(s,y) \mapsto \Esp[|Z_{s,y}|^n]^\nu$ is enough to obtain the bound $L_\cR$. 
    Thus, \ref{pr:time-Lipschitz_pinning_point} follows from inequality \eqref{eq:wasser-1_bound} alongside the fact that
    \begin{align*}
        \int_{\R} |z||\nu_{s_1,y}(z) - \nu_{s_2,y}(z)|\,\rmd z 
        &\leq \int_{\R} |z|\int_{s_1}^{s_2}|\partial_s\nu_{s,y}(z)|\,\rmd s\,\rmd z 
        \leq \int_{s_1}^{s_2} \int_{\R}|z||\partial_s\nu_{s,y}(z)|\,\rmd z\,\rmd s.
    \end{align*}
\end{proof}

\section{The equivalent optimal stopping problems}\label{sec:OSP}

Consider now the OSP
\begin{align}\label{eq:OSP}
    V(t, x) = \sup_{\tau \leq T-t}\Esp[X_{t+\tau}]_{t, x}^\nu,
\end{align}
where the supremum is taken among all stopping times of $(X_{t+u})_{u\in[0, T-t]}$, and $\Esp_{t, x}^\nu$ is the mean operator with respect to the probability measure $\Pr_{t, x}^\nu$, defined such that $\Pr[\cdot]_{t, x}^\nu = \Pr[\ \cdot\ |\ X_t = x]^\nu$. For later use, we also introduce the probability measure $\Pr_{t, x}$ such that $\Pr[\cdot]_{t, x} = \Pr[\ \cdot\ |\ X_t = x]$.

\subsection{The equivalent optimal stopping problem of a randomized Brownian bridge}

Define the transformed OSP
\begin{align}\label{eq:equivalent_OSP}
    W(s, y) = \sup_{\sigma \leq 1-s}\Esp[G(s+\sigma, Y_{s+\sigma})]_{s, y}^\nu,
\end{align}
for the gain function $G$ defined at \eqref{eq:gain_function}, and where the supremum is taken among stopping times of the process $(Y_{s+u})_{u\in[0,1-s]}$, which is a $\Pr_{s,y}$-BM starting at $Y_s = y$.
We claim that solving \eqref{eq:OSP} is equivalent to solving \eqref{eq:equivalent_OSP}, in the sense specified by the following proposition.

\begin{proposition}[Equivalence of the OSPs]\label{pr:OSP_equivalence}\ \\
    Let $V$ and $W$ be the value functions in \eqref{eq:OSP} and \eqref{eq:equivalent_OSP}. For $(t, x)\in[0,T]\times\R$, take $s = h(t)/\sT$ and $y = \frac{x - m(t)}{\varphi(t)\sqrt{\sT}}$. Then,
    \begin{align}\label{eq:equivalence_value}
        V(t, x) = W\lrp{s, y}.
    \end{align}
    Moreover, $\tau^* = \tau^*(t, x)$ is an OST for $V$ if and only if $\sigma^* = \sigma^*(s, y)$, defined such that $s + \sigma^* = h(t + \tau^*)/\sT$, is an OST for $W$.
\end{proposition}

\begin{proof} 
    From \eqref{eq:BM_representation}, we obtain the following representation for $X_{t+u}$ under $\Pr_{t, x}$:
    \begin{align*}
        X_{t+u} &= m(t+u) + \varphi(t+u)\lrp{B_{h(t+u)} + x_0} \\
        &= m(t+u) + \varphi(t+u)\sqrt{\sT}\lrp{(B_{h(t+u)} - B_{h(t)})/\sqrt{\sT} + y},
    \end{align*}
    where we used, also based on \eqref{eq:BM_representation}, that $B_{h(t)} = \sqrt{\sT}y - x_0$ under $\Pr_{t, x}$.
    Taking $r = (h(t+u) - h(t))/\sT$ and $Y_{s+r}\defeq (B_{(s+r)\sT}-B_s\sT)/\sqrt{\sT} + Y_0$, we obtain that $X_{t+u} = G\lrp{s + r, Y_{s+r}}$ under $\Pr_{s,y}^\nu$ (equivalently, under $\Pr_{t,x}^\nu$). 
    
    For every stopping time $\tau$ of $(X_{t + u})_{u\in[0, T-t]}$, define the stopping time $\sigma$ of $(Y_{s + r})_{r\in[0, 1-s]}$ such that $s + \sigma = h(t + \tau)/\sT$. 
    Hence, \eqref{eq:equivalence_value} is a consequence of
    \begin{align*}
        V(t, x) &= \sup_{\tau \leq T-t} \Esp[X_{t + \tau}]_{t, x}^\nu 
        = \sup_{\sigma \leq 1-s} \Esp[G\lrp{s + \sigma, Y_{s + \sigma}}]_{s, y}^\nu 
        = W\lrp{s, y}.
    \end{align*}
    Furthermore, suppose that $\tau^* = \tau^*(t, x)$ is an OST for \eqref{eq:OSP} and that there exists a stopping time $\sigma' = \sigma'(s, y)$ that performs better than $\sigma^* = \sigma^*(s, y)$ in \eqref{eq:equivalent_OSP}. Consider $\tau' = \tau'(t, x)$ such that $t + \tau' = h^{-1}(s + \sigma')\sT$. Then,
    \begin{align*}
        \Esp[X_{t + \tau'}]_{t, x}^\nu = \Esp[G(s + \sigma', Y_{s + \sigma'})]_{s, y}^\nu > \Esp[G(s + \sigma^*, Y_{s + \sigma*})]_{s, y}^\nu = \Esp[X_{t + \tau^*}]_{t, x}^\nu,
    \end{align*}
    which conflicts with the optimality of $\tau^*$. Analogous arguments ensure that the optimality of $\sigma^*$ in \eqref{eq:equivalent_OSP} implies that of $\tau^*$ in \eqref{eq:OSP}.
\end{proof}

\section{Analysis of the transformed optimal stopping problem}\label{sec:OSP-analysis}

Keeping in mind that the original OSP of a rGMB \eqref{eq:OSP} can be seen through the lens of an OSP of the simpler rBB \eqref{eq:equivalent_OSP}, we now obtain results on the later that immediately translate to the former. 

Denote by $\sigma^* = \sigma^*(s, y)$ to the first time the process $((s+u, Y_{s+u}))_{u\in[0,1-s]}$, starting at $(s, y)$, hits the closed (stopping) set $\cD \defeq \lrc{(s, y) : W(s,y) = G(s, y)}$. That is, 
\begin{align*}
    \sigma^* = \sigma^*(s, y) \defeq \inf\lrc{u\geq 0 : (s+u, Y_{u}^{s,y}) \in \cD},
\end{align*}
where we have defined $Y^{s,y} = \lrp{Y_{u}^{s,y}}_{u\in[0, 1-s]}$ such that 
\begin{align*}
   \mathrm{Law}\lrp{\lrp{Y_{u}^{s,y}}_{u\in[0, 1-s]}, \Pr^\nu} = \mathrm{Law}\lrp{\lrp{Y_{s+u}}_{s\in[0, 1-s]}, \Pr_{s,y}^\nu}.
\end{align*}

\begin{proposition}[Optimal stopping time characterization] \label{pr:OST_characterization}\ \\ 
    The first hitting time $\sigma^*$ is optimal in \eqref{eq:equivalent_OSP}. That is,
    \begin{align}\label{eq:OST_equivalent_OSP}
        W(s, y) = \Esp[G(s+\sigma^*, Y_{s+\sigma^*})]_{s, y}^\nu.
    \end{align}    
\end{proposition}

\begin{proof}
    On the set $\{Z_{s,y} = z\}$, the process $Y_{s+u}$ becomes a classic BB from $Y_s = y$ to $Y_1 = Z_{s,y} = z$. Then, using \eqref{eq:sup_BB} from Lemma \ref{lm:BB_results},
    \begin{align}\label{eq:bounded_sup_Y}
        \Esp[\sup_{u\leq 1-s}|Y_{s+u}|]_{s,y}^\nu 
        &= \Esp[\Esp[\sup_{u\leq 1-s}|Y_{s+u}|\ \Big|\ Z_{s,y}]_{s,y}^\nu]_{s,y}^\nu
        \leq \sqrt{2\pi}\ln(2) + |y| + \Esp[|Z_{s, y}|]^\nu < \infty,
    \end{align}
    where we recall that the random pinning point $Z_{s, y}$ has finite second moment due to Assumption \ref{asm:terminal_density}. 
    From~\eqref{eq:gain_function} along with the continuity of $a_0(s)$ and $a_1(s)$ on $[0, 1]$, we get that
    \begin{align}\label{eq:bounded_sup}
        \Esp[\sup_{u\leq 1-s}|G(s+u, Y_{s+u})|]_{s,y}^\nu < \infty,
    \end{align}
    for all $(s, y)\in[0, 1)\times\R$. 
    This regularity along with the continuity of the gain function $G$ (see, e.g., Corollary 2.9 and Remark 2.10 from \cite{Peskir-2006-optimal}) concludes the proof.
\end{proof}
 
We next explore the continuity of the value function. To this end, we need to embed two rBBs starting at different initial conditions under the same probability space. We illustrate the construction of this embedding by fixing the time component and varying the space component. Analogous arguments hold for the general case in which both components are different.  

First, consider the processes $Y^{s,y_i} = (Y_u^{s,y_i})_{u\in[0, 1-s]}$ as the same BM starting at the different initial conditions $Y_0^{s,y_i} = y_i$, $i = 1, 2$, and living in the filtered space $(\what{\Omega}, \what{\cF}, \what{\F}, \Pr)$. Introduce the random (pinning-point) variables $Z_{s,y_i}$ on the probability space $(\wtilde{\Omega}, \wtilde{\cF}, \wtilde{\Pr}^\pi)$, where $\pi$ is an admissible distribution of $(Z_{s,y_1},Z_{s,y_2})$ whose marginals are $\nu_{s,y_1}$ and $\nu_{s,y_2}$, and $\wtilde{\Pr}^\pi$ is the probability measure associated to $\pi$. Notice that we do not impose any particular form for $\pi$, thus it may be chosen to suit convenience.

Let $\Omega = \what{\Omega}\times\wtilde{\Omega}$, $\cF = \what{\cF}\otimes\wtilde{\cF}$, $\F = \F\vee\sigma(Z_{s,y_1},Z_{s,y_2})$, and $\Pr^\pi(\cdot) = \what{\Pr}\otimes\wtilde{\Pr}^\pi(\cdot \mid Y_{1-s}^{s,y_1} = Z_{s,y_1}, Y_{1-s}^{s,y_2} = Z_{s,y_2})$. Now, the process $Y^{s,y_i}$ is a rBB with random pinning point $Z_{s,y_i}$ when considered under the filtered space $(\Omega, \cF,\F,\Pr^\pi)$. Moreover, when considered jointly in the same filtered space, both processes share the same driving BM and their pinning points are jointly distributed according to $\pi$. Finally, we use $\Esp^\pi$ to refer to the mean operator with respect to $\Pr^\pi$.

Based on this embedding, we now obtain the Lipschitz continuity of the value function away from the horizon.

\begin{proposition}[Lipschitz continuity of the value function]\label{pr:value_continuous}\ \\
    For all compact sets $\cR \subset [0, 1)\times \R$, there exists a constant $L_\cR > 0$ such that
    \begin{align}\label{eq:W_Lipschitz}
       |W(s_1,y_1) - W(s_2,y_2)| \leq L_\cR(|s_1-s_2| + |y_1 - y_2|).
    \end{align}
\end{proposition}

\begin{proof}
    For an arbitrary $\varepsilon \in(0,1)$, take $s \in [0, 1-\varepsilon]$, and choose $y_1, y_2\in\R$. Denote by $\pi^* = \pi_{(s,y_1),(s,y_2)}^*$ to the Wasserstein-$1$ copula of $\nu_{s,y_1}$ and $\nu_{s,y_2}$ (see Proposition \eqref{eq:wasser-1-distance}).
    Since $|\sup_\sigma a_\sigma - \sup_\sigma b_\sigma|\leq \sup_\sigma|a_\sigma - b_\sigma|$, and due to Jensen's inequality,
    \begin{align}
        |W(s, y_1) - W(s, y_2)| 
        &= \lrav{\sup_{\sigma\leq 1-s}\Esp[G(s+\sigma, Y_{s+\sigma})]_{s,y_2}^\nu 
        - \sup_{\sigma\leq 1-s}\Esp[G(s+\sigma, Y_{s+\sigma})]_{s,y_2}^\nu} \nonumber \\
        &\leq \Esp[\sup_{u\leq 1-s} \lrav{G(s+u, Y_{u}^{s, y_1}) - G(s+u, Y_{u}^{s, y_2})}]^{\pi^*} \nonumber \\
        &\leq A_1\Esp[\sup_{u\leq 1-s}\lrav{Y_{u}^{s,y_1} - Y_{u}^{s,y_2}}]^{\pi^*}, \label{eq:value_diff_space}
    \end{align}
    where we recall the shape of $G$ in \eqref{eq:gain_function} and the bounding constants in \eqref{eq:bounds} to obtain the last inequality. 
    We can use \eqref{eq:sup_BBs_diff_space} from Lemma \ref{lm:BB_results}, as well as Proposition \ref{pr:wasser-distance}, to get the following bound: 
    \begin{align}
        \Esp[\sup_{u\leq1-s}\lrav{Y_{u}^{s,y_1} - Y_{u}^{s,y_2}}]^{\pi^*}
        &= \Esp[\Esp[\sup_{u\leq1-s}\lrav{Y_{u}^{s,y_1} - Y_{u}^{s,y_2}}\ \mid \ Z_{s,y_1}, Z_{s,y_2}]^{\pi^*}]^{\pi^*} \nonumber \\
        &\leq |y_1-y_2| + \Esp[|Z_{s,y_1} - Z_{s,y_2}|]^{\pi^*} \leq L_{\varepsilon}^{(1)}|y_1-y_2|,
        \label{eq:sup_rBBs_diff_space}
    \end{align}
    where $L_\cR^{(1)}$ is a positive constant.

    Take now $(s_1, y)$ and $(s_2, y)$ in a compact set $\wtilde{\cR}\subset[0, 1)\times \R$, with $s_1 \leq s_2$. Reasoning as above, we obtain the following for $\pi^*~=~\pi_{(s_1,y),(s_2,y)}^*$,
    \begin{align}
        |W(s_1, y) &- W(s_2, y)| \nonumber \\
        &= \lrav{\sup_{\sigma\leq 1-s_1}\Esp[G(s_1+\sigma, Y_{s_1+\sigma})]_{s_1,y}^\nu 
        - \sup_{\sigma\leq 1-s_2}\Esp[G(s_2+\sigma, Y_{s_2+\sigma})]_{s_2,y}^\nu} \nonumber \\
        &= \lrav{\sup_{\sigma\leq 1-s_1}\Esp[G(s_1+\sigma, Y_{s_1+\sigma})]_{s_1,y}^\nu 
        - \sup_{\sigma\leq 1-s_1}\Esp[G(s_2+\sigma\wedge(1-s_2), Y_{s_2+\sigma\wedge(1-s_2)})]_{s_2,y}^\nu} \nonumber \\
        &\leq \Esp[\sup_{u\leq 1-s_1} \lrav{G(s_1+u, Y_{u}^{s_1, y}) - G(s_2+u\wedge(1-s_2), Y_{u\wedge(1-s_2)}^{s_2, y})}]^{\pi^*} \nonumber \\
        &\leq\lrp{A_0' + A_1'\Esp[\sup_{u\leq 1-s_1}|Y_u^{s_1,y}|]^\nu}(s_2-s_1) + A_1\Esp[\sup_{u\leq1-s_1}\lrav{Y_{u}^{s_1,y} - Y_{u\wedge(1-s_2)}^{s_2,y}}]^{\pi^*}. \label{eq:value_diff_time}
    \end{align}
    Hence, by relying on \eqref{eq:sup_BBs_diff_time} from Lemma \ref{lm:BB_results}, recalling the bound \eqref{eq:bounded_sup_Y} and the uniform boundedness of $\Var[Z_{s,y}]$ in Assumption \ref{asm:terminal_density}, and using Proposition \ref{pr:wasser-distance}, we obtain that
    \begin{align}
        \Esp[\sup_{u\leq1-s_1}\lrav{Y_{u}^{s_1,y} - Y_{u\wedge(1-s_2)}^{s_2,y}}]^\nu 
        &= \Esp[\Esp[\sup_{u\leq1-s_1}\lrav{Y_{u}^{s_1,y} - Y_{u\wedge(1-s_2)}^{s_2,y}}\ \Big |\ Z_{s_1,y}, Z_{s_2,y}]^{\pi^*}]^{\pi^*} \nonumber \\
        &\leq (s_2-s_1)L_{\wtilde{\cR}}^{(2)} + \Esp[|Z_{s_1,y} - Z_{s_2,y}|]^{\pi^*} 
        \leq L_{\wtilde{\cR}}^{(2)}(s_2-s_1), \label{eq:sup_rBBs_diff_time}
    \end{align}
    for some positive constants $L_{\wtilde{\cR}}^{(2)}$ and $L_{\wtilde{\cR}}^{(3)}$.
    Hence, the proof of \eqref{eq:W_Lipschitz} concludes after plugging \eqref{eq:sup_rBBs_diff_space} and \eqref{eq:sup_rBBs_diff_time} back into \eqref{eq:value_diff_space} and \eqref{eq:value_diff_time}, respectively. 
\end{proof}

Denote by $\InfGen_Y$ the infinitesimal generator of the time-space process $((s,Y_s))_{s\in[0, 1]}$ under $\Pr^\nu$. That is, for a suitable function $f$,
\begin{align}\label{eq:InfGen}
    (\InfGen_Y f)(s, y) = \partial_s f(s,y) + \mu^\nu(s, y)\partial_y f(s, y) + \frac{1}{2}\partial_{yy} f(s, y).
\end{align}
Introduce the (continuation) open set $\cC \defeq \cD^c = \{(s, y) : W(s, y) > G(s, y) \}$.

The so-far proven regularity of the value function, along with the smoothness of the drift and diffusion coefficients in \eqref{eq:SDE_rBB}, suffice to prove that the solution of the OSP \eqref{eq:equivalent_OSP} also solves a Free-Boundary Problem (FBP) with $\InfGen_Y$ as the (parabolic) differential operator and $\partial\cD$ as the free boundary, called the Optimal Stopping Boundary (OSB) in optimal stopping lingo. Further details about this connection between OSPs and FBPs can be found in \citet[Chapter III, Section 7]{Peskir-2006-optimal}. 

\begin{proposition}[The free-boundary problem]\ \\
    The value function satisfies $W\in C^{1, 2}(\cC)$ and, together with the OSB $\partial \cD$, it solves the following FBP
    \begin{subequations}
    \begin{align}
        \InfGen_Y W &= 0\   \text{\quad on } \cC  \label{eq:FBP_diff_operator} \\
        W &= G   \text{\quad on } \partial \cD \quad (\text{instantaneous stop})\label{eq:FBP_instant_stopping}
    \end{align}
    \end{subequations}
\end{proposition}

\begin{proof}
    Since $W$ is continuous on $\cC$ (see Proposition \ref{pr:value_continuous}) and the coefficients of the parabolic operator $\InfGen_Y$ satisfy sufficient smoothness conditions (local $\alpha$-Hölder continuity suffice) classical results in parabolic partial differential equations \cite[Section 3, Theorem 9]{Friedman-1964-partial} guarantee that, for any open rectangle $\cR \subset \cC$, the boundary value problem
    \begin{align}\label{eq:PDE}
        \begin{cases}
            \InfGen_Y f = 0, & \text{in } \cR,  \\
            f = W, & \text{on } \partial \cR,
        \end{cases}
    \end{align}
    admits a unique solution $f \in C^{1,2}(\cR)$. 

    Applying Itô's formula to $f(s + u, Y_{s+u})$ and evaluating at $u = \sigma_{\cR} = \inf\{u\geq 0 : (s+u, Y_{s+u})\notin \cR\}$, and then taking expectations under  $\Pr_{s, y}^\nu$-expectation for $(s, y) \in \cR$ and using the strong Markov property, we finally deduce that
    $
    W(s, y) = \Esp[W(s + \sigma_{\cR}, Y_{s + \sigma_{\cR}})]_{s, y}^\nu = f(s, y).
    $
\end{proof}

\begin{remark}
    The FBP \eqref{eq:FBP_diff_operator}-\eqref{eq:FBP_instant_stopping} has only one boundary condition, which is typically called ``instantaneous stopping'' in optimal stopping terminology. To guarantee uniqueness of solution an extra boundary condition is needed, which typically, comes in the form of smoothly pasting the value and gain functions at the boundary, hence its name ``smooth-fit'' condition. Specifically, it requires $\partial_y W(s, y) = \partial_y G(s, y)$ for all $(s, y) \in \partial \cD$.

    Obtaining the smooth-fit condition is crucial in optimal stopping theory, but also mathematically challenging. As a rule of thumb, if the process enters the stopping set immediately after starting at the OSB, the smooth-fit condition holds (see, e.g., \cite{DeAngelis-2020-global} and \citet[Section 9.1]{Peskir-2006-optimal}). The common approach to obtain this (probabilistic) regularity requires smoothness of the OSB.

    Our framework, however, poses many extra challenges compared to typical settings that prevented us from obtaining the smooth-fit condition. The most notorious difficulty is the absence of characterization of $\partial \cD$ as the graph set of time-dependent functions. Even in case such characterization exists, the highly nonlinear space dependence of the drift adds another layer of complexity in applying standard techniques employed to obtain the common regularities of $\partial \cD$: monotonicity plus continuity (see e.g., \cite{Peskir-2005-American}); continuity and piecewise monotonicity (see. e.g., \citet[Example 7]{DeAngelis-2020-global} and \cite{Friedman-1964-partial}); Lipschitz continuity (see Remark 4.5 from \cite{DeAngelis-2019-Lipschitz} for a comment on this approach and Proposition 8 in \cite{Azze-2024-optimal-GMB} for an implementation).
\end{remark}

We next compare the value functions associated to different prior terminal densities ordered in the \textit{likelihood ratio} sense. 
A random variable $Z_1$ with density $\nu_1$ is said to be smaller, in the likelihood ratio order, than another (independent) random variable $Z_2$ with density $\nu_2$, 
denoted by $Z_1 \leqlr Z_2$ (or, in an abuse of notation $\nu_1 \leqlr \nu_2$), if 
\begin{align}\label{eq:likelihood_ratio_condition}
    \Pr[Z_1 \in \cA]\Pr[Z_2 \in \cB] \leq \Pr[Z_1 \in \cB]\Pr[Z_2 \in \cA],
\end{align}
for all Borel sets $\cA \leq \cB$, in the sense $z_1 \leq z_2$ whenever $z_1\in\cA$ and $z_2\in\cB$. 

\begin{proposition}[Ordering of value functions]\label{pr:value_order}
    Consider two densities $\nu_1$ and $\nu_2$ and denote by $W_1$ and $W_2$ to their associated value functions \eqref{eq:equivalent_OSP}, respectively. Then, 
    \begin{align*}
        \nu_1 \leqlr \nu_2 \Longrightarrow W_1 \leq W_2.
    \end{align*}
    Consequently, $\cD_2 \subset \cD_1$ and $\cC_1 \subset \cC_2$.
\end{proposition}

\begin{proof}
    Let $Z_{i,s,y}$ be two independent pinning points, for $i = 1,2$, associated to the posterior densities $\nu_{i,s,y}$ (given $Y_s = y$) of the priors $\nu_i$, for $i = 1, 2$, as defined in \eqref{eq:updated_terminal_dens}. Then, for two Borel sets $\cA \leq \cB$ and taking $f(z) = \frac{\phi(z;y,1-s)}{\phi(z;x_0,1)}$,
    \begin{align*}
        &\Pr[Z_{1,s,y} \in \cA]^{\nu_1}\Pr[Z_{2,s,y} \in \cB]^{\nu_2} 
        - \Pr[Z_{1,s,y} \in \cB]^{\nu_1}\Pr[Z_{2,s,y} \in \cA]^{\nu_2} \\
        &= \int_{\cA}\int_{\cB}f(z)f(z')\nu_{2,s,y}(z')\nu_{1,s,y}(z)\, \rmd z'\rmd z 
        - \int_{\cA}\int_{\cB}f(z)f(z')\nu_{1,s,y}(z')\nu_{2,s,y}(z)\rmd z'\rmd z \\
        &= \int_{\cA}\int_{\cB}f(z)f(z')\lrp{\nu_{2,s,y}(z')\nu_{1,s,y}(z) - \nu_{1,s,y}(z')\nu_{2,s,y}(z)}\rmd z'\rmd z \\
        &\leq \max_{\substack{z\in\cA \\ z'\in\cB}}\lrc{f(z)f(z')}\int_{\cA}\int_{\cB}\lrp{\nu_{2,s,y}(z')\nu_{1,s,y}(z) - \nu_{1,s,y}(z')\nu_{2,s,y}(z)}\rmd z'\rmd z \leq 0,
    \end{align*}
    which means that $Z_{1,s,y} \leqlr Z_{2,s,y}$ and, consequently (see, e.g., Theorem 1.C.1 by \cite{Shaked-2007-stochastic}), $\Pr[Z_{1,s,y} \geq z]^{\nu_1} \leq \Pr[Z_{2,s,y} \geq z]^{\nu_2}$ for all $z\in\R$, which results in $\Esp[Z_{1,s,y}]^{\nu_1} \leq \Esp[Z_{2,s,y}]^{\nu_2}$ and, therefore,
    \begin{align*}
        \mu^{\nu_1}(s, y) \leq \mu^{\nu_{2}}(s, y),
    \end{align*}
    for all $(s, y) \in[0, 1)\times\R$, with $\mu^{\nu_i}$ as in \eqref{eq:drift_rBB}. By a standard comparison argument in SDEs (see e.g., Corollary 3.1 from \cite{Peng-2006-necessary}), the corresponding rBBs $Y_{u}^{s, y,(i)}$ hold the relation $Y_{u}^{s, y,(1)}\leq Y_{u}^{s, y,(2)}$ $\Pr^\pi$-a.s. for all $u\in[0, T]$, where $\Pr^{\pi}$ can be consider as the probability measure induced by a common underlying BM and the joint density of the pinning points $Z_{i,s,y}$, which we are assuming to be independent and, hence $\pi(z_1,z_2) = \nu_{1,s,y}(z_1)\nu_{2,s,y}(z_2)$. See the construction of this probability space in the prelude of Proposition \ref{pr:value_continuous}.
    
    It follows (recall that $y \mapsto G(s, y)$ is increasing) that $W_1 \leq W_2$ and, therefore, $\cD_2 \subset \cD_1$ and $\cC_1 \subset \cC_2$. 
\end{proof}

A direct consequence of Proposition \ref{pr:value_order} is that the value function of priors with lower and/or upper bounded support can be controlled by the value function of a Gauss--Markov Bridge (GMB), which is comprehensively studied in \cite{Azze-2024-optimal-GMB} (see Section \ref{sec:dirac-delta_prior}).

\begin{corollary}[Bounds for the value function with bounded-supported priors] \label{cor:bounds_value_bounded-supported}\ \\
    Let $W_{z^*}$ be the value function associated to taking $\nu_{z^*} = \delta_{z^*}$, for $z^*\in\R$. That is, $Y$ is a BB with pinning point $z^*$ under $\Pr^{\nu_{z^*}}$. Let $\cD_{z^*}$ and $\cC_{z^*}$ be the corresponding stopping and continuation sets, respectively. Then,
    \begin{enumerate}[label=(\textit{\roman{*}}), ref=(\textit{\roman{*}})]
        \item if $\supp\,\nu \subseteq (-\infty,z^*]$, then $W \leq W_{z^*}$ and, consequently, $\cD_{z^*}  \subset \cD$, \label{cor:upper-bounded_support}
        \item if $\supp\,\nu \subseteq [z^*,\infty)$, then $W \geq W_{z^*}$ and, consequently, $\cC_{z^*} \subset \cC$.\label{cor:lower-bounded_support}
    \end{enumerate}
\end{corollary}

\begin{proof}
    Under \ref{cor:upper-bounded_support} (\textit{resp.} \ref{cor:lower-bounded_support}), the inequality relating the gain functions is a direct consequence of Proposition \ref{pr:value_order}, after checking that $\nu \leqlr \delta_{z^*}$ (\textit{resp.} $\delta_z^* \leqlr \nu$). Then,
    \begin{align*}
        (s, y) \in \cD_{z^*} &\Rightarrow W(s, y) \leq W_{z^*}(s, y) = G(s,y) \Rightarrow (s, y) \in \cD \\
        (\textit{resp. } (s, y) \in \cC_{z^*} &\Rightarrow W(s, y) \geq W_{z^*}(s, y) > G(s,y) \Rightarrow (s, y) \in \cC).
    \end{align*}
\end{proof}


Finally, in the next proposition, we derive some useful maximal bounds related to the BB that we used before in previous results. 

\begin{lemma}[Brownian bridge maximal bounds]\label{lm:BB_results}\ \\ 
    For a BB $\lrp{Y_{u}^{s,y,z}}_{u\in[0, 1-s]}$ from $Y_{0}^{s,y,z} = y$ to $Y_{1-s}^{s,y,z} = z$, the following inequalities hold: 
    \begin{align}
        \Esp[\sup_{u\in[0, 1-s]}|Y_u^{s, y, z}|] &\leq \sqrt{1-s}\sqrt{\pi/2}\ln(2) + |y| + |z|, \label{eq:sup_BB} \\
        \Esp[\sup_{u\in[0, 1-s]}|Y_u^{s, y_1, z_1} - Y_{u}^{s, y_2, z_2}|] &\leq |y_1 - y_2| + |z_1 - z_2|, \label{eq:sup_BBs_diff_space} \\ 
        \Esp[\sup_{u\in[0, 1-s_1]}|Y_u^{s_1, y, z_1} - Y_{u\wedge(1-s_2)}^{s_2, y, z_2}|] 
        &\leq |z_1-z_2| +  \frac{s_2-s_1}{(1-s_1)(1-s_2)}\lrp{\frac{|z_2|+|y| + 5\sqrt{\pi/2}}{4}}, \label{eq:sup_BBs_diff_time}   
    \end{align}
    where, in the last inequality, we take $0 \leq s_1\leq s_2 < 1$, and assume that both BBs are driven by the same BM.
\end{lemma}

\begin{proof}
    Let $\lrp{Y_r}_{r\in[0, 1]}$ be a Brownian bridge from $Y_0 = 0$ to $Y_1 = 0$. Using the well-known result (see, e.g., formulas (4.4) and (4.5) from \cite{Biane-2001-probability})
    \begin{align*}
        \Pr[\sup_{0\leq r\leq 1}|Y_r| \geq x] = 2\sum_{n=1}^\infty (-1)^{n-1}e^{-2 n^2 x^2},
    \end{align*}
    for $x>0$, we obtain, by integrating the tail probability and interchanging the sum and integration operators, that
    \begin{align*}
        \Esp[\sup_{0\leq r\leq 1}|Y_r|] 
        = 2\sum_{n=1}^\infty (-1)^{n-1}  \int_0^\infty e^{-2 n^2 x^2}\,\rmd x 
        = \sqrt{\pi/2}\sum_{n=1}^\infty (-1)^{n-1}/n 
        = \sqrt{\pi/2}\ln(2).
    \end{align*}
    Hence, \eqref{eq:sup_BB} follows after considering the decomposition
    \begin{align}\label{eq:BB_standard_representation}
        Y_u^{s,y,z} = Y_{r}/\sqrt{1-s} + y + (z-y)r,\quad r = u/(1-s).
    \end{align}

    Obtaining \eqref{eq:sup_BBs_diff_space} is trivial after
    recalling our assumption that both BBs are coupled under the same BM, and thus they share the same stochastic part in the anticipative representations $Y_{u}^{s,y_i,z_i} = y_i + (z_i - y_i)\frac{u}{1-s} + \lrp{B_u - \frac{u}{1-s}B_{1-s}}$.

    To tackle \eqref{eq:sup_BBs_diff_time} we introduce the change of time $h(s) = s/(s-1)$, take $t = h(s)$ and $r = h(s+u) - h(s)$, and rely on the following well-known (see, e.g., Equation 2 and the proof of Proposition 2 in \cite{Azze-2024-optimal-GMB}) representation of a BB,
    \begin{align*}
        Y_{u}^{s,y,z} = \frac{1}{1+t+r}\lrp{zr + B_r + y(1+t)},
    \end{align*}
    where $(B_r)_{r\in\R_+}$ is a standard BM. Hence, for $t_i = h(s_i)$, we get that
    \begin{align*}
        &\Esp[\sup_{u\in[0, 1-s_1]}|Y_u^{s_1, y, z_1} - Y_{u\wedge(1-s_2)}^{s_2, y, z_2}|] \\
        &\leq \sup_{r\in\R_+}\lrav{\frac{z_1r}{1+t_1+r} - \frac{z_2r}{1+t_2+r}} 
        + |y|\sup_{r\in\R_+}\lrav{\frac{1+t_1}{1+t_1+r} - \frac{1+t_2}{1+t_2+r}} 
        + \Esp[\sup_{r\in\R_+}\lrav{\frac{B_r}{1+t_1+r} - \frac{B_r}{1+t_2+r}}],
    \end{align*}
    where, for each addend, we can verify the following bounds: 
    \begin{align*}
       \sup_{r\in\R_+}\lrav{\frac{z_1r}{1+t_1+r} - \frac{z_2r}{1+t_2+r}} 
       & \leq |z_1-z_2|\sup_{r\in\R_+}\lrav{\frac{r}{1+t_1+r}} + |z_2|\sup_{r\in\R_+}\lrav{\frac{r}{1+t_1+r} - \frac{r}{1+t_2+r}} \\
       & \leq |z_1-z_2| + |z_2||t_1-t_2|\sup_{r\in\R_+}\frac{r}{(1+r)^2}, \\
       |y|\sup_{r\in\R_+}\lrav{\frac{1+t_1}{1+t_1+r} - \frac{1+t_2}{1+t_2+r}} 
       & \leq |y||t_1-t_2|\sup_{r\in\R_+}\frac{r}{(1+r)^2}, \\
       \Esp[\sup_{r\in\R_+}\lrav{\frac{B_r}{1+t_1+r} - \frac{B_r}{1+t_2+r}}] 
       & \leq |t_1-t_2|\sup_{r\in\R_+}\frac{|B_r|}{(1+r)^2} 
    \end{align*}
    Hence, \eqref{eq:sup_BBs_diff_time} follows after using the fact that $\sup_{r\in\R_+} r/(1+r)^2 = 1/4$, recalling that $t_2 - t_1 = \frac{s_2-s_1}{(1-s_1)(1-s_2)}$, and relying on the following inequality:
    \begin{align*}
        \Esp[\sup_{r\in\R_+}\lrav{B_r}\frac{r}{(1+r)^2}] 
        &\leq \Esp[\sup_{r\in[0, 1]}\lrav{B_r}\frac{r}{(1+r)^2}] + \Esp[\sup_{r\in(1, \infty)}\lrav{B_r}/r] \\
        &\leq \frac{1}{4}\Esp[\sup_{r\in[0, 1]}\lrav{B_r}] + \Esp[\sup_{r\in(1, \infty)}\lrav{B_{1/r}}] 
        = \frac{5}{4}\Esp[\sup_{r\in[0, 1]}\lrav{B_r}] <  \frac{5}{4}\sqrt{\frac{\pi}{2}},
    \end{align*}
    where we used the fact that $r/(1+r)^2 \leq 1/r$ for all $r > 0$ in the first inequality, the time-inversion property of the BM for the equality, and the bound $\Esp[\sup_{r\in[0, 1]}\lrav{B_r}] \leq \sqrt{\pi/2}$ for the last inequality.
\end{proof}

\section{The single boundary case}\label{sec:single_boundary}

Next we provide sufficient conditions under which $\partial \cD$ can be represented as the graph of a single function.

\begin{proposition}[Sufficient condition for single optimal stopping function]\label{pr:single-boundary_stopping} \ \\
The following propositions hold true:
   \begin{enumerate}[label=(\textit{\roman{*}}), ref=(\textit{\roman{*}})]
        \item If $y\mapsto \mu^\nu(s,y)$ is non-decreasing for all $s\in[0,1)$ or, equivalently, $\Var[Z_{s,y}] \geq (1-s)$, and if $a_1$ is non-decreasing (meaning that $\beta$ in \eqref{eq:SDE_GM} is positive), then there exists a function $b:[0, 1]\rightarrow[-\infty,\infty]$ such that $\cD = \lrc{(s, y) : y \leq b(s)}$. \label{pr:single-boundary_stopping_below}
        \item If $y\mapsto \mu^\nu(s,y)$ is non-increasing for all $s\in[0,1)$ or, equivalently, $\Var[Z_{s,y}] \leq (1-s)$, and if $a_1$ is non-increasing (meaning that $\beta$ in \eqref{eq:SDE_GM} is negative), then there exists a function $b:[0, 1]\rightarrow [-\infty,\infty]$ such that $\cD = \lrc{(s, y) : y \geq b(s)}$.\label{pr:single-boundary_stopping_above}
    \end{enumerate}
\end{proposition}

\begin{proof}
    First, a direct differentiation of \eqref{eq:drift_rBB} yields 
    \begin{align*}
        \partial_y \mu^\nu(s, y) = \frac{1}{1-s}\lrp{\frac{\Var[Z_{s,y}]}{1-s} - 1},
    \end{align*}
    which clarifies the equivalence between the monotonicity of $\mu^\nu$ with respect to $y$ and the comparison of the posterior variance against $(1-s)^{-1}$.

    To obtain \ref{pr:single-boundary_stopping_below}, 
    apply Itô's formula to $G(s+u, Y_{s+u})$, evaluate at a given stopping time $u = \sigma \leq 1-s$ and remove the martingale term after taking $P^\nu$-expectation, and then take the supremum over $\sigma$, which leaves us with the equation
    \begin{align*}
        W(s, y) - G(s, y) 
        &= \sup_{\sigma \leq 1-s}\Esp[\int_0^{\sigma} (\InfGen_Y G)(s+u, Y_{u}^{s,y})]^\nu.
    \end{align*}
    Since $y\mapsto (\InfGen_Y G)(s,y) = a_0'(s) + a_1'(s)y + \mu^\nu(s,y)a_1(s)$ is non-decreasing for all $s\in[0, 1)$, and since $Y_{u}^{s,y_1} \leq Y_{u}^{s,y_2}$ for all $u\in[0, 1-s]$ $\Pr^\nu$-a.s. whenever $y_1 \leq y_2$ (see, e.g., Theorem 1.5.5.9 from \cite{Jeanblanc-2009-mathematical}), then, $y \mapsto (W-G)(s, y)$ is also non-decreasing for all $s\in[0, 1)$. Consequently, if $(s, y_1) \in \cD$, meaning that $(W-G)(s, y_1) = 0$, then $(W-G)(s, y_2) \leq 0$ whenever $y_1\geq y_2$, meaning that $(s,y_2) \in \cD$.

    The proof of \ref{pr:single-boundary_stopping_above} follows analogous arguments.
\end{proof}

We now introduce a family of prior densities for which an application of Proposition \ref{pr:single-boundary_stopping} guarantees the existence of a single boundary separating the stopping and continuation sets: 
A density $\nu$ is said to be $\gamma^2$-strongly log-concave if there exists a log-concave function $f$ for which $\nu(z) = f(z)\phi(z,0, \gamma^2)$. We denote by $SLC_{\gamma^2}$ to the class of $\gamma^2$-strongly log-concave densities. 

\begin{corollary}[$SLC_1$ priors implies single boundaries] \label{cor:strong-log-concave} \ \\
   If $\nu \in SLC_1$ and $a_1$ is a non-increasing function (meaning that $\beta$ in \eqref{eq:SDE_GM} is non-positive), there exists a function $b:[0, 1]\rightarrow [-\infty,\infty]$ such that $\cD = \lrc{(s, y) : y \geq b(s)}$.
\end{corollary}

\begin{proof}
    If $\nu \in SLC_1$, then the function $r(z) = \nu(z)/\phi(z;x_0,1)$ is log-concave. Consequently (see, e.g., Theorem 1.1 from \cite{Harge-2004-convex}, taking $g(x) = x^2$), since $\nu_{s,y}(z) \propto r(z)\phi(z, y, 1-s)$, then $\Var[Z_{s,y}] \leq (1-s)$, which concludes the proof after applying \ref{pr:single-boundary_stopping_above} from Proposition \ref{pr:single-boundary_stopping}.
\end{proof}

Note that the proof of Corollary \ref{cor:strong-log-concave} also guarantees that the class of priors $SLC_1$ meets Assumption \ref{asm:terminal_density}, enabling all the results in Section \ref{sec:OSP-analysis}. The class $SLC_1$ can be seen as imposing more informative priors than the unconditioned BM, in the sense of choosing a unimodal prior that reduces the variance and lightens the tails of the unconditioned terminal standard normal density. We refer to \cite{Saumard-2014-log} for a exhaustive literature review on log-concave and strongly log-concave distributions. 

Next, we highlight two particular types of priors within the $SLC_1$ that have been already tackled in the optimal-stopping literature: the Gaussian and the degenerated priors.

\subsection{The Gaussian prior}

Consider now the prior terminal density $\nu$ to be normal, that is, $\nu(z) = \phi(z;\theta,\gamma^2)$ for some $\theta \in \R$ and $\gamma^2 > 0$. In such case, by identifying the ratio and multiplication of normal densities as another scaled normal density (see Lemma~\ref{lm:gaussians}), the drift function in \eqref{eq:drift_rBB} takes the form 
$\mu^\nu(s, y) = (\theta_{s, y} - y)/(1-s)$, for
\begin{align*}
    \theta_{s,y} &= \gamma_{s}^2\lrp{\frac{y}{1-s} + \frac{\theta}{\gamma^2} - x_0}, \quad 
    \gamma_{s}^2 = \lrp{\frac{1}{1-s} + \frac{1}{\gamma^2} - 1}^{-1}.
\end{align*}
That is, $\mu^\nu(s, y) = a(s) + b(s)y$, with
\begin{align*}
    a(s) &= \frac{\gamma_s^2}{1-s}\lrp{\frac{\theta}{\gamma^2} - x_0}, \quad
    b(s) = \frac{1}{1-s}\lrp{\frac{\gamma_{s}^2}{1-s} - 1}.
\end{align*}
Hence, the gain process $G_s = G(s, Y_s)$ solves the SDE
$dG_s = \lrp{A(s) + B(s)G_s}\rmd s + a_1(s)\rmd B_s^\nu$, 
for
\begin{align*}
    A(s) &= a_0'(s) - \left(\frac{a_1'(s)}{a_1(s)} + b(s)\right)a_0(s) + a_1(s)a(s), \quad 
    B(s) = \frac{a_1'(s)}{a_1(s)} + b(s).
\end{align*}

Assume that $a_1'(s) < 0$ for all $s \in [0, 1]$ (if and only if $\beta(t) < 0$ for all $t \in [0, T]$ in \eqref{eq:SDE_GM}), and that $\gamma^2 < 1$. In such case, the OSP \eqref{eq:equivalent_OSP} (and thus \eqref{eq:OSP}) has a single boundary in form specified in \ref{pr:single-boundary_stopping_above} from Proposition \ref{pr:single-boundary_stopping}. Moreover, under these conditions, $B(s) < 0$ for all $s\in[0, 1]$ and, hence the SDE of $G_s$ can be reformulated as 
\begin{align}\label{eq:SDE_rBB_gaussian-prior}
    dG_s = \wtilde{B}(s)\lrp{\wtilde{A}(s) - G_s}\rmd s + a_1(s)\rmd B_s^\nu,
\end{align}
where $\wtilde{A}(s) \defeq -A(s)/B(s)$ and $\wtilde{B}(s) \defeq -B(s)$. 

These types of processes were considered in \cite{Azze-2024-optimal-TDOU} to solve the problem of optimally exercising an American option. 
Indeed, by taking a null strike price and interest rate in the work of \cite{Azze-2024-optimal-TDOU}, and checking Remarks 1 and 2 \textit{ibid.}, one obtains the solution of the OSP $W^+(s,y) = \sup_{\sigma\leq 1-s} \Esp[(G_{s+\sigma})^+]_{s,y}$.
Although the solution associated to a positive-part payoff does not extend immediately to an identity gain, the methodology employed by \cite{Azze-2024-optimal-TDOU} still works without major changes. Actually, the analysis is simpler due to the removing of the discounting factor and of the non-differentiability of the positive-part  gain function. 
Hence, the OSB $\mathsf{b}:[0,1]\rightarrow \R$ of the OSP \eqref{eq:equivalent_OSP}, with $G_s$ as in \eqref{eq:SDE_rBB_gaussian-prior}, is the unique solution, among the class of continuous, bounded variation functions, of the Volterra-type integral equation
\begin{align}\label{eq:free-boundary_equation_rBB_gaussian-prior}
    \mathsf{b}(s) &= \Esp[G_{1}]_{s,b(s)}^\nu + \int_s^1 \Esp[\wtilde{B}(u)(\wtilde{A}(u) - G_u)\Ind(G_u \geq \mathsf{b}(u))]_{s,\mathsf{b}(s)}^\nu\,\rmd u,
\end{align}
with $\mathsf{b}(1) = \wtilde{A}(1)$.

\subsection{The Dirac delta prior}\label{sec:dirac-delta_prior}

Consider now a Dirac-delta prior terminal density, that is, $\nu = \delta_{z^*}$ for some $z^*\in\R$. In such case, $Y = (Y_s)_{s\in[0, 1]}$ is a BB with pinning point $z^*$ under $\Pr^\nu$ and, hence, the process $(G_s)_{s\in[0, 1]} = G(s, Y_s)$ satisfies
\begin{align*}
    \Esp[G_s] &= a_0(s) + a_1(s)(x_0 + (z^* - x_0)s),\\
    \Cov[G_s,G_{s'}] &= R_1(\min(s,s'))R_2(\max(s,s')), 
\end{align*}
for $R_1(s) = a_1(s)s$ and $R_2(s) = a_1(s)(1-s)$. Therefore, $(G_s)_{s\in[0, 1]}$ is a GMB in the sense specified in \cite{Azze-2024-optimal-GMB}. If we furthermore assume that $a_1$ and $a_2$ are twice continuously differentiable (which boils down to ask for continuously differentiable coefficients $\alpha$, $\beta$ and $\zeta$ in \eqref{eq:SDE_GM}), the results in \cite{Azze-2024-optimal-GMB} imply (see Section 5 \textit{ibid.}) that $\cD = \{(s, y) : y \geq \mathsf{b}(s)\}$, for the OSB $\mathsf{b}:[0, 1]\rightarrow \R$ characterized as the unique solution, among the class of continuous functions of bounded variation, of the free-boundary equation
\begin{align}\label{eq:free-boundary_equation_rBB_dirac-prior}
    \mathsf{b}(s) &= z^* - \int_s^1 \Esp[B(u)(A(u) - G_u)\Ind(Y_u \geq \mathsf{b}(u))]_{s,\mathsf{b}(s)}^\nu\,\rmd u,
\end{align}
for 
\begin{align*}
    A(s) &= \frac{\lrp{a_{0}'(s) - \frac{a_{1}'(s)}{a_{1}(s)}a_{0}(s)}(1-s) + a_{1}(s)z^* + a_{0}(s)}{B(s)(1-s)}, \quad
    B(s) = \frac{1}{1-s} - \frac{a_1'(s)}{a_1(s)},    
\end{align*}
such that $\rmd G_s = B(s)(A(s) - G_s)\rmd s + a_1(s)\rmd B_s^\nu$.

\section{Numerical results}\label{sec:numerics}

We now explore the solution of the OSP \eqref{eq:OSP} for various sets of coefficients ($\alpha$, $\beta$, and $\zeta$ in \eqref{eq:SDE_GM}) as well as different prior densities. While the auxiliary problem \eqref{eq:equivalent_OSP} has been valuable for simplifying notation and establishing the theoretical results in Sections \ref{sec:OSP-analysis} and \ref{sec:single_boundary}, it is more practical to compute the numerical solution of \eqref{eq:OSP} directly in the original $(t, x)$ coordinates. This approach avoids the extra step of solving the transformed problem and then converting the solution back via the equivalence stated in Proposition \ref{pr:OSP_equivalence}. In particular, the time change $s = s(t)$ defined in \eqref{eq:change-time} accounts for $s'(t) = \frac{\zeta^2(t)}{\varphi^2(t)}$. Thus, depending on the slope and volatility coefficients, $s'(t)$ can become extremely large or small. Such behavior may distort the time partitions in the original $t$-domain into highly uneven ones in the transformed $s$-domain, complicating the balance between the two time discretizations. Additionally, the spatial coordinate transformation $y = y(t, x)$ is not time-homogeneous, which means that rectangular grids in the $(t, x)$ space are mapped into irregular shapes in the transformed plane.

We implement Monte Carlo's simulations and use the dynamic principle to recursively identify stopping and continuation points. That is, given a time-space grid $\{(t_i, x_j) : i = 1,\dots, N, j = 1, \dots, M\}$ of the rectangle $[0, T]\times [\uline{x},\oline{x}]$, for some $N$, $M \in \N$, and for $\uline{x}$, $\oline{x} \in \R$ such that $x_{0} \in [\uline{x},\oline{x}]$, the following algorithm is implemented to return a boolean matrix $D$ whose $(i,j)$-th entry equals $1$ if $(t_i, x_j)$ is classified as a stopping point, as well as an approximation of the value function $V$. Details are laid out in Algorithm \ref{alg:stopping_continuation}. We highlight the use of $\wtilde{\nu}$ instead of $\nu$ to refer to the prior of the unconditioned GMP in \eqref{eq:SDE_GM}, that is $\Pr_{0, x_0}^{\nu}(X_T = \wtilde{Z}) = 1$, for $\wtilde{Z}\sim \wtilde{\nu}$. Due to \eqref{eq:BM_representation}, $\wtilde{\nu}(z) = \nu((z - a_0(1))/a_1(1))/a_1(s)$. 


\begin{algorithm}[H]
\caption{Optimal stopping solver}
\begin{algorithmic}[1]
\Require GMP coefficients: $\alpha$, $\beta$, $\zeta$; space boundaries $\uline{x}$, $\oline{x}$; time and space grids $\pmb{t}~=~(t_i)_{i=0}^N$ and $\pmb{x} = (x_i)_{i=1}^M$; prior $\wtilde{\nu}$; Monte Carlo sample size $K$; initial condition $x_0$.  
  \State Initialize the value and decision matrices $(V_{i,j})_{i,j=1}^{N,M}$ and $(D_{i,j})_{i,j=1}^{N,M}$.
  \State Set $V_{N,j} = x_{j}$ and $D_{N,j} = 1$ for all $j = 1,\dots, M$.
  \For{$i = (N-1):1$}
    \For{$j = 1:M$}
      \State Sample from posterior, $z_k \sim \wtilde{\nu}_{t_i,x_{j}}$ for $k = 1,\dots,K$. 
      \State Simulate $K$ observations of a BB from $(t_i, x_{j})$ to $(T, z_k)$ at time $t_{i+1}$:
             $B_k \sim \mathcal{N}(\theta, \gamma^2)$, for 
             \begin{align*}
                \theta &= m_{t_i, x_j}(t_{i+1})  +  \Cov[X_{t_{i+1}},X_{T}] (z - m_{t_i, x_j}(T) ) / v_{t_i}(T) , \\
                \gamma^2 &= v_{t_i}(t_{i+1}) -  (\Cov[X_{t_{i+1}},X_{T}])^2 / v_{t_i}(T). 
             \end{align*}
      \State Construct estimations $\wtilde{V}_k^j \approx V(t_{i+1}, B_k)$ via interpolation of the pairs $((x_{j}, V_{i+1,j}))_{j=1}^M$.
      \State Estimate the continuation value $V_\text{next}^j = \frac{1}{K}\sum_{k=1}^K\wtilde{V}_k^j$.
      \State Set $V_{i,j} = \max\{x_j, \,V_\text{next}^j\}$.
      \State Set $D_{i,j} = 1$ if $V_{i,j} = x_{j}$ and $D_{i,j} = 0$ otherwise.
    \EndFor
  \EndFor
\State \textbf{Output:} value-function matrix $V$; decision matrix $D$.
\end{algorithmic}
\label{alg:stopping_continuation}
\end{algorithm}

All the figures displayed in this section are produced by running Algorithm \ref{alg:stopping_continuation} with the space bounds $\uline{x} = -3$ and $\oline{x} = 3$, and the time and space grids $t_i = \ln(1 + i\frac{e - 1}{N})$ and $x_j = \uline{x} + j\frac{\oline{x} - \uline{x}}{M}$, for $N = M = 10^3$. We also set a Monte Carlo sample size of $K = 10^4$.

To illustrate the complex geometry of the stopping and continuation regions even in a relatively simple setting, we display in Figure \ref{fig:two_points_bridge} the output of Algorithm \ref{alg:stopping_continuation} for $\wtilde{\nu} \equiv (\delta_{-1} + \delta_{1})/2$. That is, $\Pr[X_T = 1]_{0, x_0}^\nu = \Pr[X_T = -1]_{0,x_0}^\nu = 1/2$. Besides showing disconnected stopping/continuation regions, we highlight the strong dependency on the initial condition $x_0$, and recall that this dependency is missing in the deterministic BB case.

\ifthenelse{\boolean{figcolor}}
{%
    \begin{figure}[ht]
        \centering
        \subfloat[$x_0 = -1$.]{%
        \resizebox*{0.32\textwidth}{!}{\includegraphics{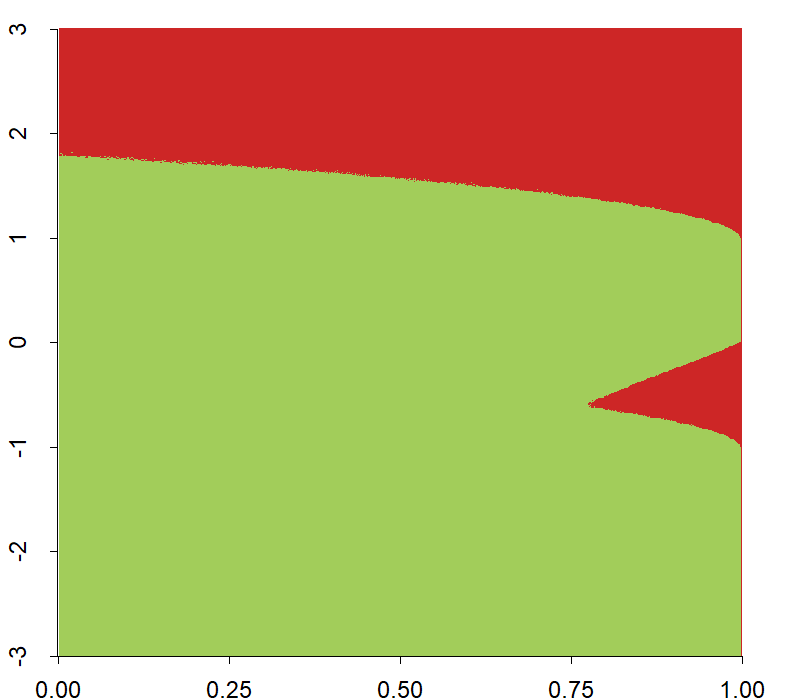}}}\hspace{2pt}
        \subfloat[$x_0 = 0$.]{%
        \resizebox*{0.32\textwidth}{!}{\includegraphics{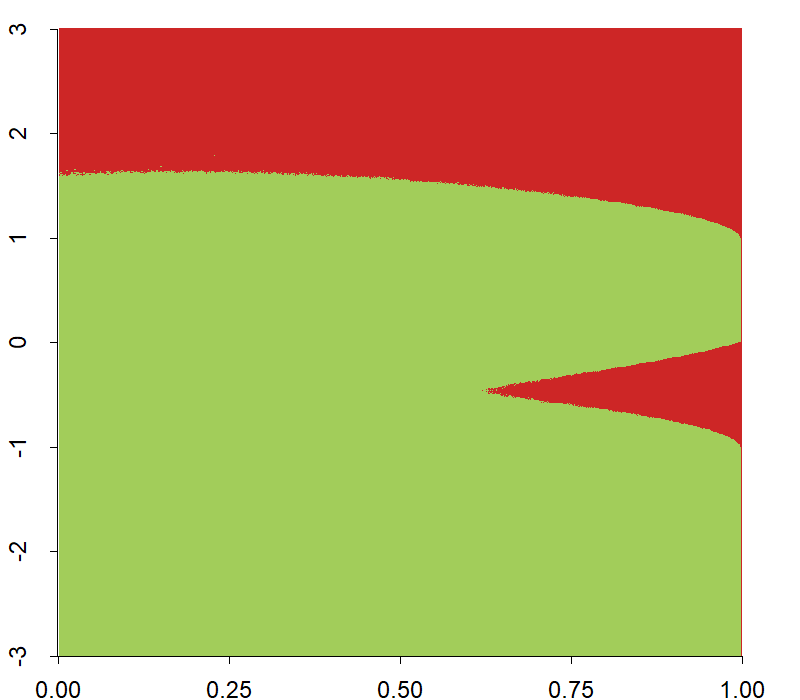}}}\hspace{2pt}
        \subfloat[$x_0 = 1$.]{%
        \resizebox*{0.32\textwidth}{!}{\includegraphics{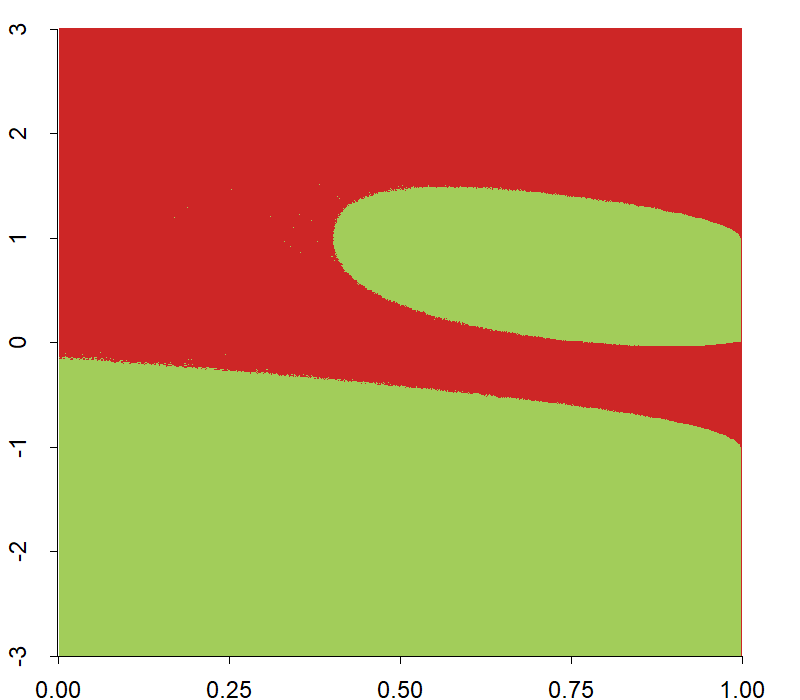}}}
        \caption{\small Numerical computation of the stopping and continuation regions for a BM forced to hit either $(1,1)$ or $(1,-1)$, each with equal probability, which corresponds to choosing the coefficients $\alpha \equiv 0$, $\beta\equiv 0$, and $\zeta \equiv 1$ in \eqref{eq:SDE_GM}, and the prior terminal density $\wtilde{\nu}(z) = (\delta_{-1} + \delta_{1})/2$. The initial condition $x_0$ is specified in the caption of each image. The red and green areas indicate the stopping and continuation sets, respectively.
        } 
        \label{fig:two_points_bridge}
    \end{figure}
}
{%
    \begin{figure}[ht]
        \centering
        \subfloat[$x_0 = -1$.]{%
        \resizebox*{0.32\textwidth}{!}{\includegraphics{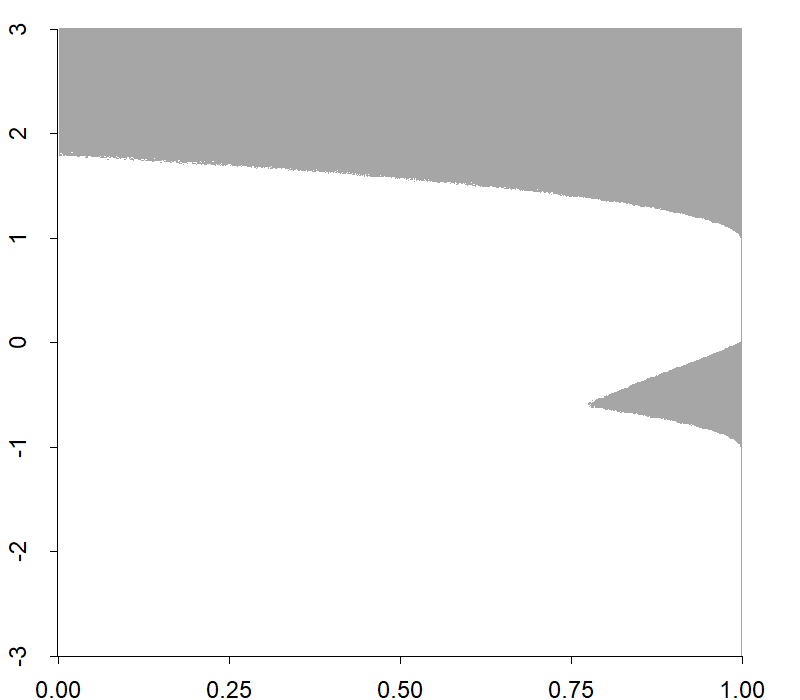}}}\hspace{2pt}
        \subfloat[$x_0 = 0$.]{%
        \resizebox*{0.32\textwidth}{!}{\includegraphics{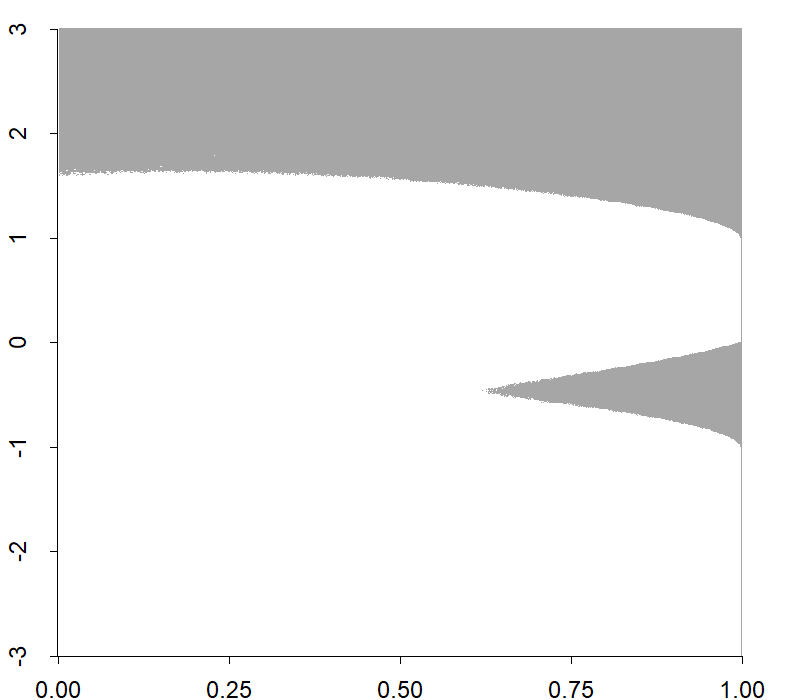}}}\hspace{2pt}
        \subfloat[$x_0 = 1$.]{%
        \resizebox*{0.32\textwidth}{!}{\includegraphics{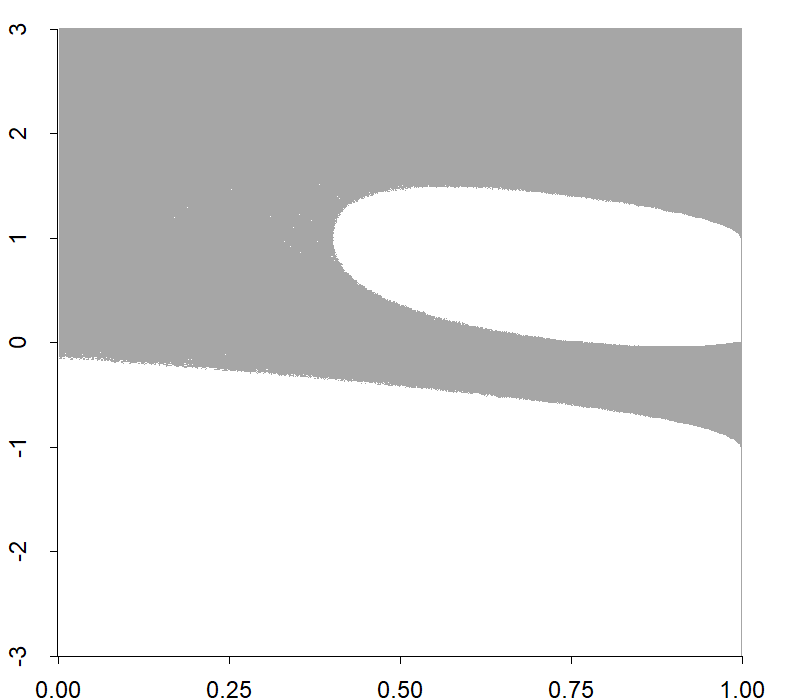}}}
        \caption{\small Numerical computation of the stopping and continuation regions for a BM forced to hit either $(1,1)$ or $(1,-1)$, each with equal probability, which corresponds to choosing the coefficients $\alpha \equiv 0$, $\beta\equiv 0$, and $\zeta \equiv 1$ in \eqref{eq:SDE_GM}, and the prior terminal density $\wtilde{\nu}(z) = (\delta_{-1} + \delta_{1})/2$. The initial condition $x_0$ is specified in the caption of each image. The gray and white areas indicate the stopping and continuation sets, respectively.
        } 
        \label{fig:two_points_bridge}
    \end{figure}
}

To provide visual insight on how the time-dependencies of the unconditioned GMP's coefficients affect the shape of the OSB, we display in Figure \ref{fig:coefficients_effect} three time-dependent parameter regimes, each motivated from a modeling point of view: image \ref{fig:coefficients_effect}-a sets $\alpha(t) = A\,\sin(\omega\pi t)$ for $A \in \R$, and $\beta \equiv -1$  to capture mean reversion towards cyclical or seasonal levels; image \ref{fig:coefficients_effect}-b models a pulling force towards the level $\alpha/\beta\equiv 0$ that rapidly (but smoothly) transitions from low ($\theta_0$) to high ($\theta_1$), by setting $\beta(t) = \theta_0 + \frac{\theta_1 - \theta_0}{2}(1 + \tanh(k(t - t_0)))$; finally, image \ref{fig:coefficients_effect}-c portraits a deterministic volatility smile regime by considering $\zeta(t) = \gamma_0 + p(t - t_0)^k$. For all of these examples we impose a zero-mean, $v_0(T)/2$-variance normal prior density, whose support is truncated in the positive semi-axis. The variance $v_0(T)/2$ ensures the resulting truncated normal prior density is a super-informative prior in the $SLC_{v_0(T)}$ sense stated in Corollary \ref{cor:strong-log-concave}, but recalling that we are working now in original coordinates $(t, x)$. Accordingly, since $\beta$ remains negative across all the cases, the OSBs are the graph of a single function.

\ifthenelse{\boolean{figcolor}}
{%
    \begin{figure}[ht]
        \centering
        \subfloat[
        $\beta(t) \equiv -1,\quad \zeta(t)\equiv 1$
        ]{%
        \resizebox*{0.32\textwidth}{!}{\includegraphics{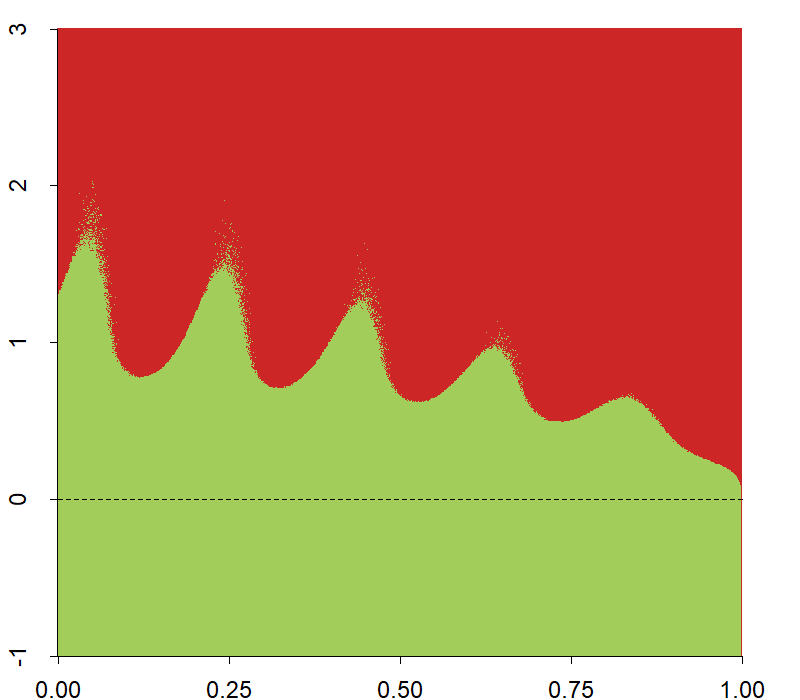}}}\hspace{2pt}
        \subfloat[
        $\alpha(t) \equiv 0,\quad \zeta(t)\equiv 1$
        ]{%
        \resizebox*{0.32\textwidth}{!}{\includegraphics{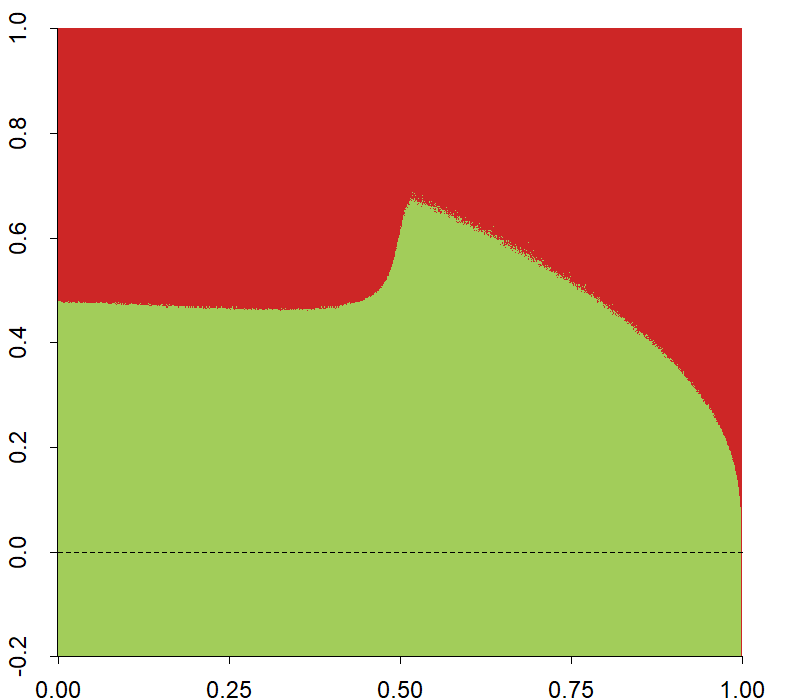}}}\hspace{2pt}
        \subfloat[
        $\alpha(t) \equiv 0,\quad \beta(t)\equiv -1$
        ]{%
        \resizebox*{0.32\textwidth}{!}{\includegraphics{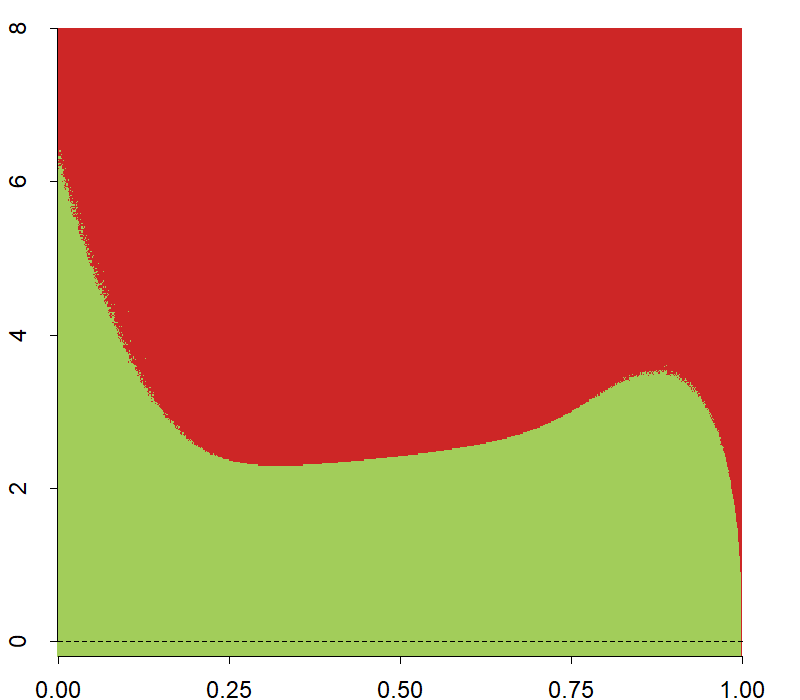}}}
        \caption{\small Numerical computation of the stopping and continuation regions associated to a GMP conditioned to adopt a positively-truncated normal terminal density $\nu(z) \propto \Ind(z\geq 0)\phi(z;0,v_0(T)/2)$. Each image highlights the time-dependent effect of one coefficient, while keeping the other constant: $\alpha(t) = 2\sin(10\pi t)$ is used in image (a); image (b) considers $\beta(t) = -10 + 0.475(1 + \tanh(100(t - 0.5)))$; and image (c) takes $\zeta(t) = 0.25 + (4(t-0.5))^4$. $x_0 = 0$ was used for all images. The red and green areas indicates the stopping and continuation sets, respectively.  
        } 
        \label{fig:coefficients_effect}
    \end{figure}
}
{%
    \begin{figure}[ht]
        \centering
        \subfloat[
        $\beta(t) \equiv -1,\quad \zeta(t)\equiv 1$
        ]{%
        \resizebox*{0.32\textwidth}{!}{\includegraphics{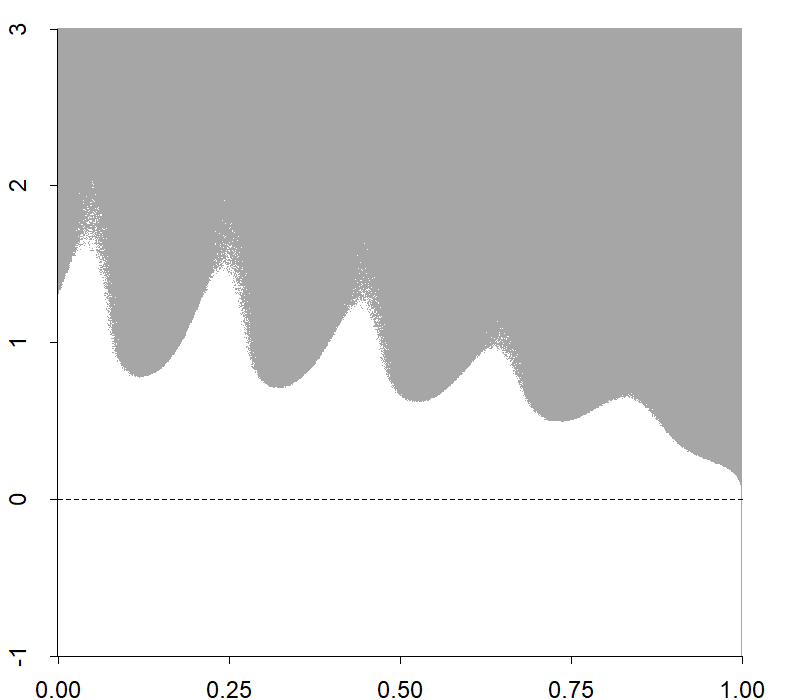}}}\hspace{2pt}
        \subfloat[
        $\alpha(t) \equiv 0,\quad \zeta(t)\equiv 1$
        ]{%
        \resizebox*{0.32\textwidth}{!}{\includegraphics{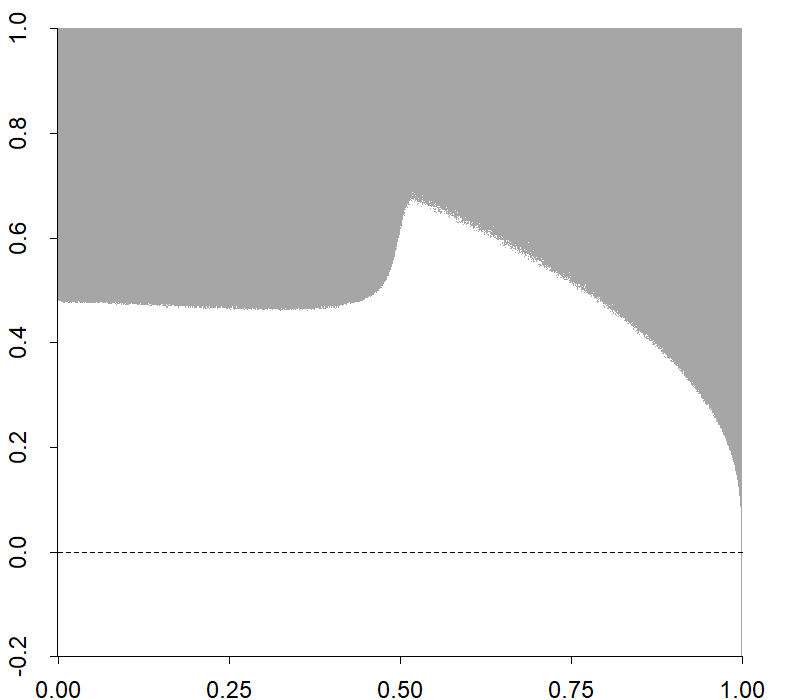}}}\hspace{2pt}
        \subfloat[
        $\alpha(t) \equiv 0,\quad \beta(t)\equiv -1$
        ]{%
        \resizebox*{0.32\textwidth}{!}{\includegraphics{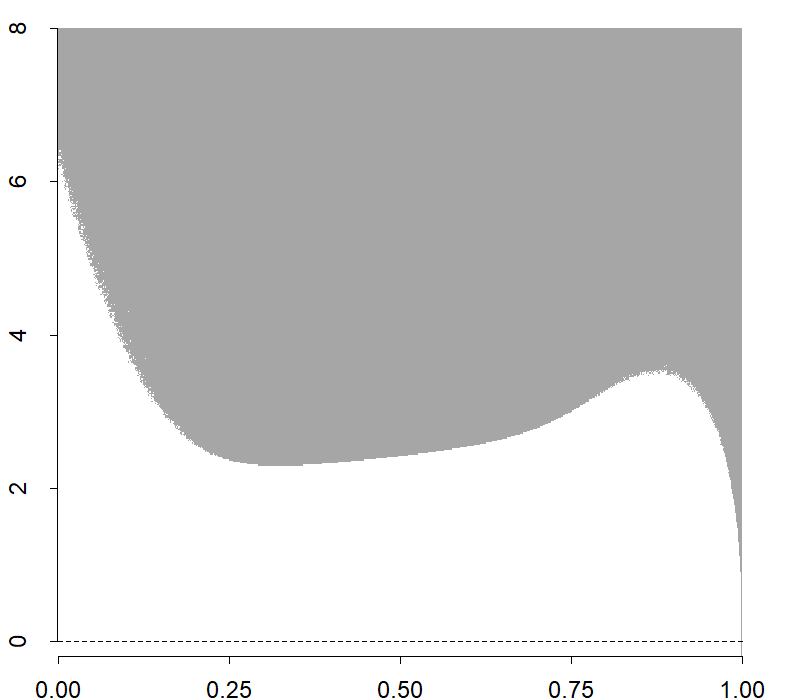}}}
        \caption{\small Numerical computation of the stopping and continuation regions associated to a GMP conditioned to adopt a positively-truncated normal terminal density $\nu(z) \propto \Ind(z\geq 0)\phi(z;0,v_0(T)/2)$. Each image highlights the time-dependent effect of one coefficient, while keeping the other constant: $\alpha(t) = 2\sin(10\pi t)$ is used in image (a); image (b) considers $\beta(t) = -10 + 0.475(1 + \tanh(100(t - 0.5)))$; and image (c) takes $\zeta(t) = 0.25 + (4(t-0.5))^4$. $x_0 = 0$ was used for all images. The gray and white areas indicates the stopping and continuation sets, respectively.  
        } 
        \label{fig:coefficients_effect}
    \end{figure}    
}

All the code necessary to implement Algorithm \ref{alg:stopping_continuation} and reproduce Figures \ref{fig:two_points_bridge} and \ref{fig:coefficients_effect} are provided in the GitHub repository \url{https://github.com/aguazz/OSP_rGMB_MC}.

\section*{Appendix}

\begin{proposition}\label{pr:bridge_SDE}
    Let $Y = (Y_s)_{s\in[0, 1]}$ be a standard BM. Then, under the probability $\Pr^\nu$ defined in Section \ref{sec:terminal_density}, $Y$ is the unique strong solution of the SDE \eqref{eq:SDE_rBB}.
\end{proposition}

\begin{proof} 
    Denote by $p_{s_0, y_0}(s, \cdot)$ and $p_{s_0, y_0}^\nu(s, \cdot)$ to the transition densities of $Y$ under $\Pr$ and $\Pr^\nu$, respectively, starting at $(s_0, y_0)$ and transitioning to $(s, \cdot)$. 
    
    The Markovianity of $Y$ implies that
    \begin{align*}
        p_{0, x_0}^{\nu}(s, y) &= \int_\R p_{0, x_0}^{z}(s, y)\nu(z)\, \rmd z 
        =  p_{0, x_0}(s, y)\int_\R\frac{p_{s, y}(1, z)}{p_{0,x_0}(1, z)}\nu(z)\, \rmd z 
        = p_{0, x_0}(s, y)\psi^{\nu}(s, y),
    \end{align*} 
    where $\psi^{\nu}(s, y)$ is defined at \eqref{eq:RN_derivative} and $p_{0, x_0}^{z}(s, \cdot)$ is the transition density of a BB with pinning point $z$.
    
    This means that $(\psi^{\nu}(s, Y_s))_{s\in(0,1)}$ is the Radon–Nikodym derivative process $\rmd \Pr^{\nu} / \rmd \Pr$.
    Hence (see, e.g., \citet[Theorem 3.4, Chapter III]{Jacod-2007-limit}), since $\Pr^\nu$ is absolute continuous with respect to $\Pr$ (their Radon–Nikodym derivative is positive), the density process 
    $(\psi^{\nu}(s, Y_s))_{s\in(0,1)}$ is a $\Pr$-martingale. 
    Therefore, the drift term in the Itô expansion of $\psi^{\nu}(s, Y_s)$ vanishes, resulting in
    \begin{align*}
        \rmd \psi^{\nu}(s, Y_s) 
        &= \partial_y \psi^{\nu}(s, Y_s) \rmd B_s 
        = \frac{\partial_y \psi^{\nu}(s, Y_s)}{\psi^{\nu}(s, Y_s)} \psi_{0,y}^{\nu}(s, Y_s) \rmd B_s \\
        &= \partial_y \ln\lrp{\psi^{\nu}}(s, Y_s)\psi^{\nu}(s, Y_s) \rmd B_s.
    \end{align*}
    Hence, Girsanov's theorem yields \eqref{eq:SDE_rBB}.
\end{proof}

\begin{lemma}\label{lm:gaussians}
    Let $f_i$ be a normal density with mean $\theta_i$ and variance $\gamma_i^2$, for $i = 1, 2, 3$, such that $\gamma_1 < \gamma_3$. Consider the constants
    \begin{align*}
        \gamma^2 = \lrp{\frac{1}{\gamma_1^2} + \frac{1}{\gamma_2^2} - \frac{1}{\gamma_3^2}}^{-1}, \quad
        \theta = \gamma^2\lrp{\frac{\theta_1}{\gamma_1^2} + \frac{\theta_2}{\gamma_2^2} - \frac{\theta_3}{\gamma_3^2}}, \quad
        C = \frac{\theta_1^2}{2\gamma_1^2} + \frac{\theta_2^2}{2\gamma_2^2} - \frac{\theta_3^2}{2\gamma_3^2} - \frac{\theta^2}{2\gamma^2}.
    \end{align*}
    Then,
    \begin{align}\label{eq:f1_times_f2_over_f3}
        f_1(z)f_2(z)/f_3(z) = \frac{\gamma_3}{\gamma_1 \gamma_2} \gamma e^{-C}f(z)
    \end{align}
    where $f$ is a normal density with mean $\theta$ and variance $\gamma^2$. 
\end{lemma}

\begin{proof}
    Notice that
    \begin{align*}
        f_1(z)f_2(z)/f_3(z)
        &= \frac{\gamma_3}{\sqrt{2\pi}\gamma_1\gamma_2}\exp\lrc{- \lrp{\frac{(z-\theta_1)^2}{2\gamma_1^2} + \frac{(z-\theta_2)^2}{2\gamma_2^2} - \frac{(z-\theta_3)^2}{2\gamma_3^2}}}.
    \end{align*}
    Also,
    \begin{align*}
        \frac{(z-\theta_1)^2}{2\gamma_1^2} + \frac{(z-\theta_2)^2}{2\gamma_2^2} - \frac{(z-\theta_3)^2}{2\gamma_3^2} 
        &= \frac{z^2}{2}\lrp{\frac{1}{\gamma_1^2} + \frac{1}{\gamma_2^2} - \frac{1}{\gamma_3^2}} - 
        z\lrp{\frac{\theta_1}{\gamma_1^2} + \frac{\theta_2}{\gamma_2^2} - \frac{\theta_3}{\gamma_3^2}} + 
        \lrp{\frac{\theta_1^2}{2\gamma_1^2} + \frac{\theta_2^2}{2\gamma_2^2} - \frac{\theta_3^2}{2\gamma_3^2}} \\
        &= \frac{(z - \theta)^2}{2\gamma^2} + C,
    \end{align*}
    from where it follows \eqref{eq:f1_times_f2_over_f3}.
\end{proof}

\bibliographystyle{apalike}
\bibliography{bib.bib}

\end{document}